\pdfoutput=1
\documentclass[11pt]{article}

\usepackage{amsmath,amssymb,amsthm}
\usepackage[margin=1.1in]{geometry}
\usepackage{booktabs}
\usepackage{array}
\usepackage{graphicx}
\usepackage[colorlinks=true,linkcolor=blue,citecolor=blue]{hyperref}

\usepackage{subcaption}
\usepackage{placeins}

\newtheorem{theorem}{Theorem}[section]
\newtheorem{proposition}[theorem]{Proposition}
\newtheorem{lemma}[theorem]{Lemma}
\newtheorem{corollary}[theorem]{Corollary}
\newtheorem{hypothesis}[theorem]{Hypothesis}

\theoremstyle{definition}
\newtheorem{definition}[theorem]{Definition}
\theoremstyle{remark}
\newtheorem{remark}[theorem]{Remark}
\numberwithin{equation}{section}

\newcommand{\Kc}{\mathcal{K}}
\newcommand{\Sm}[1]{\mathfrak{S}_{#1}}
\newcommand{\Kmod}{\mathfrak{K}}
\newcommand{\Pislow}{\Pi_{\Sm{\pm1}}}
\newcommand{\Piker}{\Pi_{\Kmod}}
\newcommand{\mx}{m_x}
\newcommand{\mz}{m_z}
\newcommand{\chixz}{\chi_{xz}}
\newcommand{\Rop}{\mathcal{R}}
\newcommand{\Sop}{\mathcal{S}}
\newcommand{\bQ}{\bar Q}
\newcommand{\bP}{\bar P}

\title{A coupling normal form for reversible center--center Hamiltonians
with application to the restricted three-body problem}

\author{Mingpei Lin\thanks{Corresponding author: lin.mingpei.d2@tohoku.ac.jp}\\[0.5ex]
\normalsize Advanced Institute for Materials Research, Tohoku University\\
\normalsize 2-1-1 Katahira, Aoba-ku, Sendai 980-8577, Japan}

\date{}

\begin{document}
\maketitle

\begin{abstract}
Near a reversible center--center equilibrium, the Birkhoff normal form does not describe the families bifurcating from the continued tori, which resonant normal forms capture through a reduced system.
With a reflection symmetry, we develop a coupling normal form in which the bifurcating families also appear as tori of action normal forms.
Along a selected coupling direction, the first harmonic in the angle difference is absorbed into an action-dependent counterterm.
We prove that the coupled recursion determines unique counterterm coefficients, a formal action normal form, and an exact counterterm factorization.
The factorization gives a recovery constraint, whose trivial factor corresponds to the Birkhoff normal form.
At every finite order, a symmetric pair of bifurcating branches leaves the trivial branch at a critical curve, carrying approximately invariant tori and circles with residual estimates.
In the spatial circular restricted three-body problem, these give approximate halo orbits and quasi-halo tori.
\end{abstract}

% =====================================================================
\section{Introduction}\label{sec:intro}
% =====================================================================

Near the collinear libration points of the spatial circular restricted three-body problem (CRTBP), the linearized dynamics has a hyperbolic part and two elliptic modes, planar and vertical, with frequencies close to the $1{:}1$ resonance \cite{Szebehely,MHO}.
Their nonlinear interaction organizes the Lyapunov, Lissajous, halo, and quasi-halo families, which have been described by center manifold reductions and normal forms \cite{GJMS,JM,Jorba99}, by numerical computations \cite{GomezMondelo,HaroMondelo}, and by Lindstedt--Poincar\'e series since Richardson's third-order halo approximation \cite{Richardson80}.

The starting point of this paper is the coupling-induced bifurcation mechanism of Lin and Chiba \cite{LinChiba,LLC}, which attributes the halo and quasi-halo families to the nonlinear coupling of the two modes.
A coupling coefficient $\check\eta$ enters their Lindstedt--Poincar\'e series solution, and a coupling correction term $\check\eta\,\check\Delta x$, with $\check\Delta$ a series in the amplitudes, is added to the vertical equation of motion.
The condition $\check\eta\,\check\Delta=0$ recovers the original equations to the computed order, its trivial root gives the Lissajous orbits, and the nonzero roots of the bifurcation equation $\check\Delta=0$ give the quasi-halo and halo orbits.
Such series diverge in general \cite{JM,Poincare93}, so they locate the orbit families and their bifurcation but cannot prove that these orbits and tori exist.

For the two-dimensional quasi-halo tori, existence is a small-divisor problem for the interaction of the two modes.
The normal-form treatment of such interactions has its Hamiltonian foundation in Birkhoff normalization \cite{Birkhoff27}.
Gustavson systematized resonant normalizations near elliptic equilibria \cite{Gustavson66}, the Lie-transform algorithms of Hori and Deprit made high-order computations standard in celestial mechanics \cite{Hori66,Deprit69}, and general Hamiltonian normal forms and averaging are treated in \cite{AKN,SVM,Gio}.
Parameter-dependent and quasi-periodic normal-form theories further require an explicit parameter domain and divisor conditions \cite{BHT,BHS}.
A posteriori KAM theory proves the existence of invariant tori near approximately invariant tori with explicit residual bounds \cite{DLLGJV,HaroBook}, and quasi-periodic bifurcation theory follows families of tori through bifurcations \cite{BHT,Hanssmann}.

These tools apply most directly to the invariant tori of an action normal form, that is, a normal form depending only on the actions.
Whether the bifurcating families appear as such tori depends on how the normal form treats the slow harmonics.
These terms depend on the angles of the two modes only through their difference and carry the exchange between the modes.
Nonresonant normalization removes them through divisors proportional to the frequency difference \cite{GJMS,JM} and gives the Birkhoff normal form, which does not describe the bifurcating families.
Resonant normalization retains them \cite{Hanssmann,HanssmannHoveijn,PucaccoMarchesiello} and in the CRTBP represents the halo orbits as equilibria of a reduced system \cite{CPS,CCP,Pucacco19,Luo}.

Our aim is a normal form in which the bifurcating families also appear as tori of an action normal form, so that the tools above can be applied to them.
We obtain it by giving the coupling-induced bifurcation mechanism a Hamiltonian formulation.
In this formulation, we absorb the first slow harmonic along one selected coupling direction into a Hamiltonian counterterm, the counterpart of $\check\eta\,\check\Delta x$.
In the terminology of normal form styles \cite{Murdock}, absorption is a third treatment besides removal and retention.
As in the series solutions of \cite{LinChiba}, a term is added and then required to vanish, an idea analogous to Moser's modifying terms theorem and its parameterized versions \cite{Moser67,BHS,Wagener}.
The result is a new type of normal form, the coupling normal form, for general reversible center--center Hamiltonians with a reflection symmetry.

Three difficulties shape the construction.
First, a correction of the vertical equation alone is not Hamiltonian, so we consider the auxiliary family $H+\eta\Kc$ with a coupling parameter $\eta$ and a counterterm $\Kc=\Delta W$ along the coupling direction $W$, the reversible quadratic first slow harmonic, which in the CRTBP is the first slow projection of the product $xz$.
Second, $\Delta$ must be a function of the normal-form actions expressed in reference coordinates, which depend on the normalizing map still to be constructed.
The counterterm and the map are therefore constructed by a coupled recursion that is triangular in the phase-space degree, and each absorption step reduces to one scalar equation per new coefficient, solved without division by the frequency difference.
Third, the construction must preserve both symmetries.
Symmetry propagation makes the first slow obstruction odd in $\eta$, so the quotient by $\eta$ gives unique even counterterm coefficients.

The construction yields an action normal form and an exact counterterm factorization (Theorem~\ref{thm:main}), hence the recovery constraint
\[
 \eta\,\Delta(I_1,I_2;\eta)=0,
\]
the Hamiltonian counterpart of $\check\eta\,\check\Delta=0$.
Its factor $\eta=0$ is the trivial branch, described by the Birkhoff normal form of $H$ and formed in the CRTBP by the planar Lyapunov and Lissajous families, and $\Delta=0$ is the bifurcation equation.
At every finite order a symmetric pair of roots $\pm\eta$ branches off the trivial branch at a critical curve under a nondegeneracy condition (Theorem~\ref{thm:line}), and on these branches $\eta$ is a signed branch coordinate.
At order four the critical curve is a line meeting the axis $I_2=0$ at an action proportional to the frequency difference, hence at small amplitude in the CRTBP.

On the constraint set the auxiliary Hamiltonian term vanishes, but its vector field adds a tangential correction to the internal dynamics.
At finite order the constraint yields approximately invariant tori and circles with explicit residual estimates, as required by the a posteriori arguments above.
The analytic estimates hold at a fixed truncation order on nonresonant parameter intervals, and the persistence of these objects, and hence a dynamical bifurcation, is not proved here.

Sections~\ref{sec:abstract} and~\ref{sec:linear} introduce the setting and the auxiliary quadratic family, Section~\ref{sec:style} constructs the coupling normal form, and Section~\ref{sec:reduced} develops the finite-order dynamics and the bifurcation equation.
Section~\ref{sec:crtbp} applies the theory to the spatial CRTBP, including a Sun--Earth $L_1$ example, and Section~\ref{sec:discussion} discusses the relation to \cite{LinChiba,LLC} and the open problems.

% =====================================================================
\section{Reversible center--center Hamiltonians}\label{sec:abstract}
% =====================================================================

This section fixes the setting of the general theory, a two-degree-of-freedom Hamiltonian with two distinct elliptic modes.
It introduces the reference coordinates, the frequency assumptions, and the reversing and reflection symmetries.

Let the \emph{reference coordinates} $(\bQ_1,\bP_1,\bQ_2,\bP_2)$ be real canonical coordinates on a neighborhood of the origin in $\mathbb R^4$, with symplectic form
\[
d\bQ_1\wedge d\bP_1+d\bQ_2\wedge d\bP_2 .
\]
Set $\bQ:=(\bQ_1,\bQ_2)$ and $\bP:=(\bP_1,\bP_2)$.
For differentiable functions $F$ and $G$, the Poisson bracket and the Hamiltonian vector field are
\[
\{F,G\}:=\frac{\partial F}{\partial\bQ}\frac{\partial G}{\partial\bP}
        -\frac{\partial F}{\partial\bP}\frac{\partial G}{\partial\bQ},
\qquad
X_F:=\Big(\frac{\partial F}{\partial\bP},-\frac{\partial F}{\partial\bQ}\Big).
\]
Let $H$ be a real formal or real analytic Hamiltonian
\begin{equation}\label{eq:abstract-H}
 H=H_2+\sum_{n\ge3}H_n,
 \qquad
 H_2=\frac{\omega_1}{2}(\bQ_1^2+\bP_1^2)
     +\frac{\omega_2}{2}(\bQ_2^2+\bP_2^2),
\end{equation}
where every $H_n$ is homogeneous of phase-space degree $n$ and
\begin{equation}\label{eq:abstract-frequencies}
  \omega_1>\omega_2>0,
  \qquad
  \delta:=\omega_1-\omega_2>0,
  \qquad
  g:=\omega_1+\omega_2 .
\end{equation}
We call $\delta$ the \emph{frequency difference} and $g$ the \emph{frequency sum}.

The construction uses two involutions.
The first is the anti-symplectic reversor
\begin{equation}\label{eq:abstract-R}
 \mathcal R:(\bQ_1,\bP_1,\bQ_2,\bP_2)
 \longmapsto(\bQ_1,-\bP_1,\bQ_2,-\bP_2),
\end{equation}
and the second is the symplectic reflection
\begin{equation}\label{eq:abstract-S}
 \mathcal S:(\bQ_1,\bP_1,\bQ_2,\bP_2)
 \longmapsto(\bQ_1,\bP_1,-\bQ_2,-\bP_2).
\end{equation}

\begin{hypothesis}[reversibility and $\mathbb Z_2$ reflection symmetry]
\label{hyp:abstract-symmetry}
The Hamiltonian $H$ satisfies
\begin{equation}\label{eq:abstract-symmetry}
  H\circ\mathcal R=H,
  \qquad
  H\circ\mathcal S=H .
\end{equation}
\end{hypothesis}

% =====================================================================
\section{Auxiliary quadratic family and symplectic diagonalization}\label{sec:linear}
% =====================================================================

This section constructs the quadratic part of the coupling normal form as the starting point for the recursion of Section~\ref{sec:style}.
We select the coupling direction, introduce the auxiliary quadratic family along it, and diagonalize this family by a symplectic rotation.

\subsection{Coupling direction and auxiliary quadratic family}
\label{ssec:coupling-family}

The quadratic part $H_2$ is a sum of two harmonic oscillators.
Its action--angle variables in reference coordinates are
\begin{equation}\label{eq:reference-actions}
 \bar I_j:=\frac12(\bQ_j^2+\bP_j^2),
 \qquad
 \bQ_j=\sqrt{2\bar I_j}\,\cos\bar\varphi_j,
 \qquad
 \bP_j=-\sqrt{2\bar I_j}\,\sin\bar\varphi_j,
 \qquad j=1,2.
\end{equation}
Each angle $\bar\varphi_j$ is defined where $\bar I_j>0$.
In these variables $H_2=\omega_1\bar I_1+\omega_2\bar I_2$, 
and \(\bar\varphi_j\) evolves with frequency \(\omega_j\).
Among the quadratic polynomials bilinear in the two modes, those involving only the angle difference \(\bar\varphi_1-\bar\varphi_2\) are spanned by
\[
\begin{aligned}
 \bQ_1\bQ_2+\bP_1\bP_2&=2\sqrt{\bar I_1\bar I_2}\,\cos(\bar\varphi_1-\bar\varphi_2),\\
 \bQ_1\bP_2-\bP_1\bQ_2&=2\sqrt{\bar I_1\bar I_2}\,\sin(\bar\varphi_1-\bar\varphi_2).
\end{aligned}
\]
Along the flow of $H_2$ they oscillate with the frequency difference $\delta$, and we call them the \emph{quadratic first slow harmonics}.
The reversor $\mathcal R$ sends each $\bar\varphi_j$ to $-\bar\varphi_j$, so it fixes the cosine harmonic and negates the sine harmonic.
The reflection $\mathcal S$ sends $\bar\varphi_2$ to $\bar\varphi_2+\pi$ and negates both.
Neither harmonic can therefore occur in the $\mathcal S$-even Hamiltonian $H$.
We take the $\mathcal R$-even cosine harmonic as the coupling direction, so that adding it to $H_2$ preserves reversibility.

\begin{definition}[coupling direction]\label{def:abstract-W}
The coupling direction is
\begin{equation}\label{eq:abstract-W}
 W:=\bQ_1\bQ_2+\bP_1\bP_2
 =2\sqrt{\bar I_1\bar I_2}\,\cos(\bar\varphi_1-\bar\varphi_2).
\end{equation}
It is $\mathcal R$-even and $\mathcal S$-odd.
\end{definition}

Fix a \emph{coupling scale} $\kappa>0$ and set
\begin{equation}\label{eq:abstract-d00}
 d_{00}:=-\kappa\delta.
\end{equation}
For a real \emph{coupling parameter} $\eta$, define the \emph{auxiliary quadratic family}
\begin{equation}\label{eq:H2eta}
 H_2(\eta):=H_2+\eta d_{00}W.
\end{equation}
The term $\eta d_{00}W$ is the auxiliary quadratic term, and $H_2(\eta)$ is the quadratic part of the auxiliary family introduced in Section~\ref{ssec:family}.

The parity of $W$ gives
\begin{equation}\label{eq:quadratic-equivariance}
	H_2(\eta)\circ\mathcal R=H_2(\eta),
	\qquad H_2(-\eta)=H_2(\eta)\circ\mathcal S.
\end{equation}

\begin{remark}[status of the coupling parameter]\label{rem:eta-status}
The parameter $\eta$ belongs to the auxiliary family, not to $H$, and is not a physical control parameter of the original system.
The auxiliary quadratic term $\eta d_{00}W$ couples the two modes at linear order, and $\eta$ measures the strength of this linear coupling.
At $\eta=0$ this term vanishes, but the modes remain coupled through the nonlinear terms $H_n$ with $n\ge3$.
The numerical value of $\eta$ depends on the choice of the coupling scale $\kappa$, and the rescaling
\[
 (\kappa,\eta)\mapsto(c\kappa,\eta/c),\qquad c>0,
\]
leaves the term $\eta d_{00}W$ unchanged.
\end{remark}

\subsection{Symplectic diagonalization}\label{ssec:diagonalization}
For an angle \(\theta\), let \(R_\theta\) map real canonical coordinates \((Q_1,P_1,Q_2,P_2)\) to the reference coordinates of Section~\ref{sec:abstract} by
\begin{equation}\label{eq:Rtheta-definition}
	\binom{\bQ_1}{\bQ_2}
	=\begin{pmatrix}
		\cos\theta&-\sin\theta\\
		\sin\theta& \cos\theta
	\end{pmatrix}
	\binom{Q_1}{Q_2},
	\qquad
	\binom{\bP_1}{\bP_2}
	=\begin{pmatrix}
		\cos\theta&-\sin\theta\\
		\sin\theta& \cos\theta
	\end{pmatrix}
	\binom{P_1}{P_2}.
\end{equation}
Since the same orthogonal rotation acts on the positions and on the momenta, \(R_\theta\) is symplectic. When \(\theta\) is chosen as the principal rotation angle \(\theta(\eta)\) determined below, we refer to \((Q_1,P_1,Q_2,P_2)\) as the \emph{diagonal canonical coordinates}.
Define the associated actions and angles by
\begin{equation}\label{eq:diagonal-actions}
	I_j:=\frac12(Q_j^2+P_j^2),
	\qquad
	Q_j=\sqrt{2I_j}\cos\varphi_j,
	\qquad
	P_j=-\sqrt{2I_j}\sin\varphi_j,
	\quad j=1,2.
\end{equation}
Set \(I:=(I_1,I_2)\) and \(\varphi:=(\varphi_1,\varphi_2)\), and denote the angle difference by \(\psi:=\varphi_1-\varphi_2\).

\begin{proposition}[symplectic diagonalization of the auxiliary quadratic family]\label{thm:linear}
For every $\eta\in\mathbb R$, let $\theta(\eta)\in(-\pi/4,\pi/4)$ be the principal rotation angle determined by
\begin{equation}\label{eq:rotation-angle}
 \tan(2\theta(\eta))=\frac{2\eta d_{00}}{\delta}
               =-2\kappa\eta ,
\end{equation}
and set
\begin{equation}\label{eq:sigma}
 \sigma(\eta):=\cos(2\theta(\eta))
 =\frac{1}{\sqrt{1+4\kappa^2\eta^2}} .
\end{equation}
Then $R_{\theta(\eta)}$ diagonalizes $H_2(\eta)$ in the sense that
\begin{equation}\label{eq:H2etadiag}
 H_2(\eta)\circ R_{\theta(\eta)}
 =\Omega\cdot I
 :=\Omega_1(\eta)I_1+\Omega_2(\eta)I_2 ,
 \qquad
 \Omega:=(\Omega_1,\Omega_2),
\end{equation}
with the auxiliary quadratic frequencies
\begin{equation}\label{eq:Omegat}
 \Omega_1(\eta)=\frac12\Big(g+\frac{\delta}{\sigma(\eta)}\Big),
 \qquad
 \Omega_2(\eta)=\frac12\Big(g-\frac{\delta}{\sigma(\eta)}\Big).
\end{equation}
\end{proposition}

\begin{proof}
The same symmetric matrix
\[
M=\begin{pmatrix}\omega_1& \eta d_{00} \\ \eta d_{00} &\omega_2\end{pmatrix}
\]
is the coefficient matrix of $H_2(\eta)$ in the positions $\bQ$ and in the momenta $\bP$.
The simultaneous rotation \eqref{eq:Rtheta-definition} diagonalizes $M$ when $\tan(2\theta(\eta))=2\eta d_{00}/(\omega_1-\omega_2)$, which is \eqref{eq:rotation-angle}.
Since $|\theta(\eta)|<\pi/4$, the value $\cos(2\theta(\eta))$ is positive, which gives \eqref{eq:sigma}.
The diagonal entries of the rotated matrix are $\tfrac12\big(g\pm\delta/\cos(2\theta(\eta))\big)$, which gives \eqref{eq:Omegat}.
\end{proof}

Because $d_{00}=-\kappa\delta$, the principal rotation angle $\theta(\eta)$ depends on $\eta$ only through $\kappa\eta$ and is therefore unchanged under the rescaling of Remark~\ref{rem:eta-status}.
For fixed $\kappa$, it is in one-to-one correspondence with $\eta$.

The principal rotation angle $\theta(\eta)$ is odd in $\eta$, while $\sigma(\eta)$ and the auxiliary quadratic frequencies are even.
For every angle $\theta$,
\begin{equation}\label{eq:rotation-equivariance}
 R_\theta\circ\mathcal R=\mathcal R\circ R_\theta,
 \qquad
 \mathcal S\circ R_\theta=R_{-\theta}\circ\mathcal S,
\end{equation}
and since $\theta(-\eta)=-\theta(\eta)$, the second identity gives $R_{\theta(-\eta)}\circ\mathcal S=\mathcal S\circ R_{\theta(\eta)}$.
Together with \eqref{eq:quadratic-equivariance}, these identities transport the symmetries of the auxiliary quadratic family to the diagonal canonical coordinates.
There $\mathcal R$ acts by $(Q,P)\mapsto(Q,-P)$ for every $\eta$, and $\mathcal S$ acts by $(Q_1,P_1,Q_2,P_2)\mapsto(Q_1,P_1,-Q_2,-P_2)$ from the coordinates at parameter $\eta$ to those at $-\eta$.

\begin{corollary}[frequency relations]\label{prop:freq}
For every real $\eta$,
\begin{align}
 \Omega_1+\Omega_2&=\omega_1+\omega_2=g,
 \label{eq:trace}\\
 \Omega_1-\Omega_2
 &=\frac{\delta}{\sigma(\eta)}
 =\delta\sqrt{1+4\kappa^2\eta^2}.
 \label{eq:deltaprime}
\end{align}
Consequently $\Omega_1(\eta)\ge\omega_1>\omega_2\ge\Omega_2(\eta)$, with equality only at $\eta=0$.
\end{corollary}

\begin{proof}
The identities and the inequalities follow from \eqref{eq:sigma} and \eqref{eq:Omegat}.
\end{proof}

\begin{corollary}[first slow projection of the coupling direction]
\label{prop:projection-nondegeneracy}
The transformed coupling direction is
\begin{equation}\label{eq:Wtransformed}
 W_\eta:=W\circ R_{\theta(\eta)}
 =2\sigma(\eta)\sqrt{I_1I_2}\cos\psi
 +\sin(2\theta(\eta))(I_1-I_2).
\end{equation}
Its first slow component is $\sigma(\eta)(Q_1Q_2+P_1P_2)=2\sigma(\eta)\sqrt{I_1I_2}\cos\psi$.
The coefficient $\sigma(\eta)$ is positive and even in $\eta$, whereas the coefficient $\sin(2\theta(\eta))$ of the action component is odd.
\end{corollary}

\begin{proof}
Substitution of \eqref{eq:Rtheta-definition} into \eqref{eq:abstract-W} gives
\[
 W\circ R_\theta=\cos(2\theta)\,(Q_1Q_2+P_1P_2)+\sin(2\theta)\,(I_1-I_2),
\]
and \eqref{eq:diagonal-actions} gives $Q_1Q_2+P_1P_2=2\sqrt{I_1I_2}\cos\psi$.
The parity follows from that of $\theta(\eta)$.
\end{proof}

% =====================================================================
\section{Recursive construction of the coupling normal form}\label{sec:style}
% =====================================================================

This section extends the linear normalization of Section~\ref{sec:linear} to higher phase-space degrees by a coupled recursion for the counterterm and the normalizing map.
We define the normal form style, show that the recursion is triangular, solve each step under the two symmetries, and obtain the formal coupling normal form.

\subsection{Coupling counterterm and normalizing map}\label{ssec:family}

For each $\eta$, we seek a symplectic normalizing map $\Psi$ from normal-form coordinates to reference coordinates.
Its linear part is the rotation $R_{\theta(\eta)}$ of Proposition~\ref{thm:linear}.
The normal-form coordinates keep the names $(Q,P)$ of the diagonal canonical coordinates, and their actions $I_j$ and angles $\varphi_j$ are defined by \eqref{eq:diagonal-actions}.
Define the functions on the reference phase space
\begin{equation}\label{eq:dressedactions}
 J_j:=I_j\circ\Psi^{-1},
 \qquad\text{equivalently}\qquad
 J_j\circ\Psi=I_j,
 \qquad j=1,2,
\end{equation}
that is, the normal-form actions $I_j$ expressed in reference coordinates.

Let $d_{ab}(\eta)$, $a,b\in\mathbb N_0$ with $a+b\ge1$, denote the counterterm coefficients to be determined recursively.
Define
\begin{equation}\label{eq:Deltaseries}
 \Delta(J_1,J_2;\eta)
 :=d_{00}+\sum_{a+b\ge1}d_{ab}(\eta)J_1^aJ_2^b,
 \qquad d_{00}=-\kappa\delta ,
\end{equation}
and from it the counterterm and the auxiliary family
\begin{equation}\label{eq:Heta}
 \Kc:=\Delta(J_1,J_2;\eta)W,
 \qquad
 H(\eta):=H+\eta\Kc .
\end{equation}
By \eqref{eq:dressedactions}, $\Delta(J_1,J_2;\eta)\circ\Psi=\Delta(I_1,I_2;\eta)$ depends only on the normal-form actions.
An expansion in the reference actions $\bar I_j$ of \eqref{eq:reference-actions} would not have this property, because $\bar I_j\circ\Psi$ depends on the angles.

Each $J_j=I_j\circ\Psi^{-1}$ starts at phase-space degree two, so $J_1^aJ_2^bW$ starts at degree $2(a+b)+2\ge4$ when $a+b\ge1$.
The quadratic part of $\Kc$ is $d_{00}W$, so the quadratic part of $H(\eta)$ is the auxiliary quadratic family $H_2(\eta)$ of \eqref{eq:H2eta}.
We write $H(\eta)=H_2(\eta)+\sum_{n\ge3}H_n(\eta)$ with $H_n(\eta)$ homogeneous of phase-space degree $n$.
By Proposition~\ref{thm:linear}, the rotation $R_{\theta(\eta)}$ diagonalizes $H_2(\eta)$, so
\[
 H(\eta)\circ R_{\theta(\eta)}
 =\Omega\cdot I+\sum_{n\ge3}H_n(\eta)\circ R_{\theta(\eta)}.
\]
The symbols $H$, $H_2$, and $H_n$ without $(\eta)$ always denote the original Hamiltonian and its homogeneous components.
For $\Psi$, $J_j$, $\Kc$, and the other $\eta$-dependent objects constructed below, we omit $(\eta)$ unless the dependence on $\eta$ has to be shown.

Since $\Kc$ depends on $\Psi$ through $J_j$ while $\Psi$ is produced by normalizing $H(\eta)$, the counterterm and the normalizing map must be constructed simultaneously.
Fortunately, this dependence is triangular in the phase-space degree, so both can be determined degree by degree below.

\subsection{Monomial classes and normal form style}\label{ssec:modules}

To classify the terms of functions in normal-form coordinates, we use the complex canonical coordinates
\begin{equation}\label{eq:complex-convention}
 q_j:=\frac{Q_j-iP_j}{\sqrt2}
     =\sqrt{I_j}\,e^{i\varphi_j},
 \qquad
 p_j:=-i\frac{Q_j+iP_j}{\sqrt2}
     =-i\sqrt{I_j}\,e^{-i\varphi_j},
 \qquad j=1,2.
\end{equation}
For a vector $x=(x_1,x_2)$ and $k=(k_1,k_2)\in\mathbb N_0^2$, we write $x^k:=x_1^{k_1}x_2^{k_2}$, with $0^0:=1$, and $|k|:=k_1+k_2$.
The monomials
\begin{equation}\label{eq:monomial}
 q^{k_q}p^{k_p},\qquad k_q,k_p\in\mathbb N_0^2,
\end{equation}
have phase-space degree $|k_q|+|k_p|$ and angular index $\nu:=k_q-k_p\in\mathbb Z^2$, and $q^{k_q}p^{k_p}$ depends on the angles only through $e^{i\langle\nu,\varphi\rangle}$.
It is an eigenvector of the homological operator $\mathcal L_0F:=\{\Omega\cdot I,F\}$ with eigenvalue
\begin{equation}\label{eq:eigenvalue}
 \langle k_p-k_q,i\Omega\rangle=-i\langle\nu,\Omega\rangle .
\end{equation}
The real factor $\langle\nu,\Omega\rangle$ is called the \emph{divisor} of $q^{k_q}p^{k_p}$.

\begin{definition}[monomial classes]\label{def:slowmod}
	The monomials \eqref{eq:monomial} are classified by their angular index $\nu$ as follows:
\begin{itemize}
	\item[(i)] \emph{Action monomials} have $\nu=(0,0)$ and span the action kernel
	\begin{equation}\label{eq:kernel-module}
		\Kmod=\operatorname{span}\{I_1^aI_2^b:a,b\in\mathbb N_0\}.
	\end{equation}
	\item[(ii)] \emph{Slow monomials} have $\nu=(m,-m)$ for some $m\in\mathbb Z\setminus\{0\}$ and span the corresponding slow module
	\[
	\Sm{m}=\operatorname{span}\{q^{k_q}p^{k_p}:\nu=(m,-m)\}.
	\]
	The modules $\Sm{\pm1}$ are called \emph{first slow}, and $\Sm{\pm m}$ with $m\ge2$ are called \emph{higher slow}.
	\item[(iii)] \emph{Fast monomials} have $\nu_1+\nu_2\ne0$.
\end{itemize}
\end{definition}

Monomials in $\Sm{m}$ have angular factor $e^{im\psi}$ and divisor $m\delta/\sigma(\eta)$.
Slow divisors therefore involve only the frequency difference, whereas a fast divisor also contains the frequency sum.
The terms slow and fast refer only to this frequency decomposition and do not presuppose a small frequency difference.

\begin{definition}[coupling normal form style]\label{def:style}
At every phase-space degree the monomial classes are treated as follows:
\begin{center}
\begin{tabular}{@{}llll@{}}
\toprule
class & $\nu$ & treatment & divisor \\
\midrule
action kernel $\Kmod$ & $(0,0)$ & retained & $0$ \\
first slow, $\Sm{\pm1}$ & $(\pm1,\mp1)$ & absorbed by $\Kc$ & $\pm\delta/\sigma$ \\
higher slow, $\Sm{\pm m}$ & $(\pm m,\mp m)$, $m\ge2$ & removed by $\Psi$ & $\pm m\delta/\sigma$ \\
fast & $\nu_1+\nu_2\ne0$ & removed by $\Psi$ & $\langle\nu,\Omega\rangle$ \\
\bottomrule
\end{tabular}
\end{center}
Retained monomials stay in the normal form.
The first slow component at each degree is neither retained nor removed and is called the \emph{first slow obstruction}.
It is cancelled by the counterterm $\Kc$ of \eqref{eq:Heta}, whose coefficients $d_{ab}$ are chosen for this purpose.
This treatment is called \emph{absorption}.
\end{definition}

A normal form style assigns to each degree the terms that remain in the normal form \cite{Murdock}.
Definition~\ref{def:style} adds absorption into the counterterm as a third treatment besides retention and removal.
A resonant normal form instead retains every slow module \cite{Gustavson66,Hanssmann,CPS}.

\subsection{Triangularity of the coupled recursion}\label{ssec:triangular}

To determine the counterterm and the normalizing map degree by degree, we treat all functions on phase space as series in phase-space degree, with coefficients depending on $\eta$ in one of two settings.
At the formal level, the coefficients are elements of $\mathbb R[[\eta]]$, and the even ones are those in $\mathbb R[[\eta^2]]$.
At a fixed truncation order $N\ge4$, they are real analytic on a compact interval $|\eta|\le\eta_0$ with $\eta_0>0$.
For $n\ge2$ let $\mathcal P_n$ denote the space of homogeneous polynomials of phase-space degree $n$ with such coefficients.
For a series $F$, let $\mathcal T_nF$ denote its Taylor truncation through phase-space degree $n$.

The normalizing map is built from one generator for each degree $n\ge3$.
For a generating Hamiltonian $G$, denote its time-$t$ Hamiltonian flow by $\Phi_G^t$.
Composition with its time-one map is given by the Lie series $F\circ\Phi_G^1=F+\{F,G\}+\tfrac12\{\{F,G\},G\}+\cdots$.
For generators $G_k\in\mathcal P_k$, $3\le k\le n$, define the stage-$n$ objects
\begin{equation}\label{eq:partial-objects}
 \Psi^{(2)}:=R_{\theta(\eta)},
 \qquad
 \Psi^{(n)}
 :=R_{\theta(\eta)}\circ\Phi_{G_3}^1\circ\cdots\circ\Phi_{G_n}^1,
 \qquad
 J_j^{(n)}:=I_j\circ\big(\Psi^{(n)}\big)^{-1},
\end{equation}
the stage-$n$ partial sum of \eqref{eq:Deltaseries} together with the stage-$n$ counterterm
\begin{equation}\label{eq:partial-Delta}
 \Delta^{(n)}(J_1,J_2;\eta)
 :=d_{00}+\!\!\sum_{1\le a+b\le\lfloor (n-2)/2\rfloor}\!\!
   d_{ab}(\eta)\,J_1^aJ_2^b,
 \qquad
 \Kc^{(n)}:=\Delta^{(n)}\big(J_1^{(n)},J_2^{(n)};\eta\big)\,W ,
\end{equation}
and the stage-$n$ Hamiltonian
\begin{equation}\label{eq:stage-H}
 H^{(n)}:=(H+\eta\Kc^{(n)})\circ\Psi^{(n)},\qquad n\ge2 .
\end{equation}

The next proposition shows that the generator $G_n$ and the coefficients $d_{ab}$ with $2(a+b)+2=n$ affect only terms of phase-space degree at least $n$, so that the coupled dependence is triangular in the phase-space degree.

\begin{proposition}[triangular structure of the coupled recursion]
\label{prop:triangular}
Fix $n\ge3$ and suppose that the generators $G_3,\ldots,G_n\in \mathcal P_3,\ldots,\mathcal P_n$ and the coefficients $d_{ab}$ with $2(a+b)+2\le n$ have been chosen.
Then:
\begin{itemize}
\item[(i)] $\mathcal T_{n-1}J_j^{(n)}=\mathcal T_{n-1}J_j^{(n-1)}$, and the degree-$n$ component of $J_j^{(n)}-J_j^{(n-1)}$ is $-\{I_j,G_n\}\circ R_{\theta(\eta)}^{-1}$.
Consequently $\mathcal T_mJ_j^{(n)}=\mathcal T_mJ_j^{(m)}$ for every $2\le m\le n-1$.
\item[(ii)] for every $2\le m\le n-1$, $\mathcal T_mH^{(n)}=\mathcal T_mH^{(m)}$.
\item[(iii)] the degree-$n$ component of $H^{(n)}$ equals
\[
\hat H_n+\mathcal L_0G_n +\eta\!\!\sum_{2(a+b)+2=n}\!\!d_{ab}\,I_1^aI_2^b\,W_\eta ,
\]
where the known part $\hat H_n$ is the degree-$n$ component of $H^{(n-1)}$.
\end{itemize}
\end{proposition}

\begin{proof}
(i) Since $\Psi^{(n)}=\Psi^{(n-1)}\circ\Phi_{G_n}^1$, the definition \eqref{eq:partial-objects} gives
\[
 J_j^{(n)}-J_j^{(n-1)}
 =\big(I_j\circ(\Phi_{G_n}^1)^{-1}-I_j\big)\circ\big(\Psi^{(n-1)}\big)^{-1}.
\]
By the Lie series,
\[
 I_j\circ(\Phi_{G_n}^1)^{-1}=I_j-\{I_j,G_n\}+\tfrac12\{\{I_j,G_n\},G_n\}-\cdots .
\]
The Poisson bracket of homogeneous polynomials of degrees $r$ and $s$ belongs to $\mathcal P_{r+s-2}$.
Hence $I_j\circ(\Phi_{G_n}^1)^{-1}-I_j$ starts at degree $n$, with degree-$n$ component $-\{I_j,G_n\}$.
Since $(\Psi^{(n-1)})^{-1}$ has linear part $R_{\theta(\eta)}^{-1}$, composing with it does not lower the degree of any term and turns $-\{I_j,G_n\}$ into $-\{I_j,G_n\}\circ R_{\theta(\eta)}^{-1}$.
This proves the first statement, and iterating it over the stages $m+1,\ldots,n$ proves the second.

(ii)--(iii) Likewise, the definition \eqref{eq:stage-H} gives
\[
 H^{(n)}-H^{(n-1)}
 =\big(H^{(n-1)}\circ\Phi_{G_n}^1-H^{(n-1)}\big)
 +\eta\big(\Kc^{(n)}-\Kc^{(n-1)}\big)\circ\Psi^{(n)} .
\]
We show that this difference starts at degree $n$ and determine its degree-$n$ component.
It has three contributions, from the generator $G_n$, the new counterterm terms, and the counterterm terms already present.

First, by the Lie series, $H^{(n-1)}\circ\Phi_{G_n}^1-H^{(n-1)}=\{H^{(n-1)},G_n\}+\cdots$.
Since the quadratic part of $H^{(n-1)}$ is $\Omega\cdot I$, its only term of degree at most $n$ is $\{\Omega\cdot I,G_n\}=\mathcal L_0G_n\in\mathcal P_n$.

Second, the new terms $\eta d_{ab}(J_1^{(n)})^a(J_2^{(n)})^bW$ with $2(a+b)+2=n$ occur only for even $n$.
Since $J_j^{(n)}\circ\Psi^{(n)}=I_j$ by \eqref{eq:partial-objects}, the composition of such a term with $\Psi^{(n)}$ is $\eta d_{ab}I_1^aI_2^b\,(W\circ\Psi^{(n)})$.
The quadratic part of $W\circ\Psi^{(n)}$ is $W_\eta$ of \eqref{eq:Wtransformed}, so this term starts at degree $n$ with degree-$n$ component $\eta d_{ab}I_1^aI_2^bW_\eta$.

Third, a term $\eta d_{ab}(J_1^{(n-1)})^a(J_2^{(n-1)})^bW$ already present, with $1\le a+b\le\lfloor(n-3)/2\rfloor$, changes under the update $J_j^{(n-1)}\mapsto J_j^{(n)}$.
By part (i), each difference $J_j^{(n)}-J_j^{(n-1)}$ starts at degree $n$.
The product of the $a+b$ factors $J_j$ changes by a sum of products, each containing one such difference and $a+b-1$ further factors $J_j^{(n-1)}$ or $J_j^{(n)}$ starting at degree two.
Together with the quadratic factor $W$, the change of the term therefore starts at degree $n+2(a+b)\ge n+2$.
Composition with $\Psi^{(n)}$ does not lower this degree.

Hence $H^{(n)}-H^{(n-1)}$ starts at degree $n$, so $\mathcal T_{n-1}H^{(n)}=\mathcal T_{n-1}H^{(n-1)}$.
Iterating this over the stages $m+1,\ldots,n$ gives (ii).
The degree-$n$ component of the difference is $\mathcal L_0G_n$ plus the sum of the terms $\eta d_{ab}I_1^aI_2^bW_\eta$ with $2(a+b)+2=n$.
Adding the degree-$n$ component $\hat H_n$ of $H^{(n-1)}$ gives (iii).
\end{proof}

\subsection{Homological and absorption equations}\label{ssec:homological}

Fix a degree $n\ge3$ and assume that $G_3,\ldots,G_{n-1}$ and the coefficients $d_{ab}$ with $2(a+b)+2\le n-1$ have been determined so that $\mathcal T_{n-1}H^{(n-1)}$ lies in $\Kmod$.
Let $\hat H_n$ be the known part of Proposition~\ref{prop:triangular}(iii), and let $G_n\in\mathcal P_n$ denote the unknown generator.
Let $\Piker$ and $\Pislow$ denote the projections that keep only the monomials of $\Kmod$ and of the real first slow module $\Sm{+1}\oplus\Sm{-1}$, respectively.

The style of Definition~\ref{def:style} is now written as conditions on $G_n$ and on the new coefficients $d_{ab}$ with $2(a+b)+2=n$.
In this style, $\Psi$ removes no first slow monomials, so $G_n$ contains no first slow monomials.
Since the action monomials lie in the kernel of $\mathcal L_0$ and remove nothing, we also require that $G_n$ contain none of them.
These two requirements form the \emph{normalization condition}
\begin{equation}\label{eq:gauge}
 \Piker G_n=0,
 \qquad
 \Pislow G_n=0 .
\end{equation}
Since each monomial $q^{k_q}p^{k_p}$ is an eigenvector of $\mathcal L_0$ by \eqref{eq:eigenvalue}, the term $\mathcal L_0G_n$ in Proposition~\ref{prop:triangular}(iii) then contains only fast and higher slow monomials.
By Corollary~\ref{prop:projection-nondegeneracy}, each new counterterm term splits into a first slow part and an action part,
\begin{equation}\label{eq:Wproj}
 I_1^aI_2^bW_\eta
 =2\sigma(\eta)I_1^aI_2^b\sqrt{I_1I_2}\cos\psi
 +\sin(2\theta(\eta))\,I_1^aI_2^b(I_1-I_2).
\end{equation}
Projecting the degree-$n$ component of Proposition~\ref{prop:triangular}(iii) onto the classes of Definition~\ref{def:style} therefore gives one statement for each treatment.
The fast and higher slow components vanish exactly when $G_n$ satisfies the \emph{homological equation}
\begin{equation}\label{eq:homological}
 \mathcal L_0G_n=\Piker\hat H_n+\Pislow\hat H_n-\hat H_n .
\end{equation}
The first slow component vanishes exactly when the new coefficients satisfy the \emph{absorption equation}
\begin{equation}\label{eq:absorption}
 \Pislow\hat H_n+\eta\!\!\sum_{2(a+b)+2=n}\!d_{ab}\,\Pislow(I_1^aI_2^bW_\eta)=0 .
\end{equation}
The action component is retained and equals
\[
 \Piker\hat H_n+\eta\!\!\sum_{2(a+b)+2=n}\!\!d_{ab}\sin(2\theta(\eta))\,I_1^aI_2^b(I_1-I_2).
\]
Under \eqref{eq:gauge}, the degree-$n$ component of $H^{(n)}$ thus lies in $\Kmod$ exactly when \eqref{eq:homological} and \eqref{eq:absorption} hold.

However, solving \eqref{eq:homological} requires division by the fast and higher slow divisors.
The higher slow divisors cause no difficulty, since Corollary~\ref{prop:freq} gives
\[
 \frac{|m|\delta}{\sigma(\eta)}\ge |m|\delta>0,\qquad |m|\ge2 .
\]
The fast divisors, in contrast, may vanish, and they are controlled by the following hypothesis.

\begin{hypothesis}[nonresonance on a parameter interval]\label{hyp:NR}
For the chosen $N$ and $\eta_0$, there is a constant $\gamma_N>0$ such that
\begin{equation}\label{eq:abstract-NR}
 \left|\langle\nu,\Omega(\eta)\rangle\right|
 \ge\gamma_N
\end{equation}
on $|\eta|\le\eta_0$ for the angular index $\nu$ of every fast monomial of degree at most $N$.
At the formal level, $\langle\nu,\Omega(0)\rangle\ne0$ for the angular index $\nu$ of every fast monomial, so that each fast divisor $\langle\nu,\Omega(\eta)\rangle$ is a unit of $\mathbb R[[\eta]]$.
\end{hypothesis}

An interval $|\eta|\le\eta_0$ satisfying Hypothesis~\ref{hyp:NR} for the chosen $N$ is called \emph{admissible}.

\subsection{Equivariant solvability of one recursive step}\label{ssec:solvability}

The absorption equation \eqref{eq:absorption} poses two difficulties.
For each pair $(a,b)$ the counterterm supplies a single first slow direction, while the corresponding first slow part of $\hat H_n$ has a cosine and a sine component.
Moreover, the counterterm enters with the factor $\eta$, so solving for $d_{ab}$ requires division by $\eta$.
Both difficulties are resolved by the symmetries of Hypothesis~\ref{hyp:abstract-symmetry}.
The reversor $\mathcal R$ rules out the sine component, and the reflection $\mathcal S$ makes the cosine coefficient odd in $\eta$ and hence divisible by $\eta$.
We therefore seek a recursion preserving the family symmetries
\begin{equation}\label{eq:family-equivariance}
 H(\eta)\circ\mathcal R=H(\eta),
 \qquad
 H(-\eta)=H(\eta)\circ\mathcal S.
\end{equation}
For the quadratic part $H_2(\eta)$ they reduce to \eqref{eq:quadratic-equivariance}.

By the next lemma, the real first slow module of each even degree $n$ splits into two-dimensional spaces indexed by the pairs $(a,b)$ with $2(a+b)+2=n$.
The $\Rop$-even part of each space is one-dimensional, and the counterterm supplies this direction.

\begin{lemma}[one-dimensional first slow obstruction]
\label{lem:solvability}
Fix an even degree $n\ge4$.
Then:
\begin{itemize}
\item[(i)] for each pair $(a,b)\in\mathbb N_0^2$ with $2(a+b)+2=n$, the functions $I_1^aI_2^b\sqrt{I_1I_2}\cos\psi$ and $I_1^aI_2^b\sqrt{I_1I_2}\sin\psi$ span a two-dimensional space, called the \emph{first slow space} associated with $I_1^aI_2^b$, and the degree-$n$ part of the real first slow module $\Sm{+1}\oplus\Sm{-1}$ is the direct sum of these spaces.
\item[(ii)] in each of these spaces, the $\Rop$-even part is spanned by the cosine and the $\Rop$-odd part by the sine.
\item[(iii)] $\Pislow(I_1^aI_2^bW_\eta)$ lies in the first slow space associated with $I_1^aI_2^b$ and is $2\sigma(\eta)$ times the cosine, with $\sigma(\eta)>0$.
\end{itemize}
\end{lemma}

\begin{proof}
A degree-$n$ monomial in $\Sm{+1}$ has $k_q-k_p=(1,-1)$, so $k_q=(a+1,b)$ and $k_p=(a,b+1)$ for a unique pair $(a,b)\in\mathbb N_0^2$ with $2(a+b)+2=n$.
By \eqref{eq:complex-convention} it is a constant multiple of $I_1^aI_2^b\sqrt{I_1I_2}\,e^{i\psi}$, and the monomials of $\Sm{-1}$ are, up to constant factors, the complex conjugates of these.
Taking real and imaginary parts gives (i).
Since $\Rop$ acts as $(Q,P)\mapsto(Q,-P)$, the convention \eqref{eq:diagonal-actions} shows that it sends $(\varphi_1,\varphi_2)$ to $(-\varphi_1,-\varphi_2)$ and thus $\psi$ to $-\psi$, which gives (ii).
Formula \eqref{eq:Wproj} and $\sigma(\eta)>0$ prove (iii).
\end{proof}

\begin{proposition}[equivariant solvability of one recursive step]
\label{prop:equivariant}
Fix $n\ge3$ and assume Hypotheses~\ref{hyp:abstract-symmetry} and~\ref{hyp:NR}, with $n\le N$ at finite order.
Suppose that the generators $G_3,\ldots,G_{n-1}$ and the counterterm coefficients $d_{ab}$ with $2(a+b)+2\le n-1$ satisfy the \emph{equivariance conditions}
\begin{equation}\label{eq:equivariant-class}
 G\circ\Rop=-G,
 \qquad
 G(-\eta)=G(\eta)\circ\Sop,
 \qquad
 d_{ab}(-\eta)=d_{ab}(\eta) .
\end{equation}
Then:
\begin{itemize}
\item[(i)] $\Rop\circ\Psi^{(n-1)}=\Psi^{(n-1)}\circ\Rop$ and $\Psi^{(n-1)}(-\eta)=\Sop\circ\Psi^{(n-1)}(\eta)\circ\Sop$.
\item[(ii)] the stage-$(n-1)$ Hamiltonian of \eqref{eq:stage-H} satisfies the family symmetries
\[
H^{(n-1)}\circ\Rop=H^{(n-1)}, \qquad H^{(n-1)}(-\eta)=H^{(n-1)}(\eta)\circ\Sop .
\]
In particular the known part $\hat H_n$ of Proposition~\ref{prop:triangular}(iii) satisfies both identities.
\item[(iii)] writing $\hat H_n=\sum\hat h_{k_q,k_p}(\eta)\,q^{k_q}p^{k_p}$, the coefficient $\hat h_{k_q,k_p}$ is odd in $\eta$ whenever the angular index $\nu$ of $q^{k_q}p^{k_p}$ has odd second component $\nu_2$.
In particular, if $n$ is even, then for every pair $(a,b)$ with $2(a+b)+2=n$ the coefficient $h_{ab}$ of $I_1^aI_2^b\sqrt{I_1I_2}\cos\psi$ in $\hat H_n$ admits the factorization
\begin{equation}\label{eq:hfactor}
 h_{ab}(\eta)=\eta\widetilde h_{ab}(\eta^2).
\end{equation}
\item[(iv)] if $n$ is even, then for each pair $(a,b)$ with $2(a+b)+2=n$, the absorption equation has the unique solution
\begin{equation}\label{eq:dabsolution}
 d_{ab}(\eta)
 =-\frac{\widetilde h_{ab}(\eta^2)}{2\sigma(\eta)}.
\end{equation}
Formally, $d_{ab}\in\mathbb R[[\eta^2]]$.
At finite order, $d_{ab}$ is real analytic and even on every admissible interval.
If $n$ is odd, no new counterterm coefficient occurs.
\item[(v)] the homological equation \eqref{eq:homological} has a unique solution $G_n\in\mathcal P_n$ satisfying the normalization condition \eqref{eq:gauge}, and this solution satisfies
\[
G_n\circ\Rop=-G_n, \qquad G_n(-\eta)=G_n(\eta)\circ\Sop .
\]
Together with (iv), the equivariance conditions \eqref{eq:equivariant-class} therefore hold through degree $n$.
\end{itemize}
\end{proposition}

\begin{proof}
(i) By \eqref{eq:rotation-equivariance} and the oddness of $\theta(\eta)$, the linear factor commutes with $\Rop$ and satisfies $R_{\theta(-\eta)}=\Sop\circ R_{\theta(\eta)}\circ\Sop$.

For the nonlinear factors, let $G$ be one of $G_3,\ldots,G_{n-1}$.
Since $\Rop$ is a linear anti-symplectic involution, $X_{G\circ\Rop}=-\Rop_*X_G$.
Hence $G\circ\Rop=-G$ gives $\Rop_*X_G=X_G$, and $\Rop$ commutes with $\Phi_G^1$.
Since $\Sop$ is a linear symplectic involution, $\Phi_{G(\eta)\circ\Sop}^1=\Sop\circ\Phi_{G(\eta)}^1\circ\Sop$, and $G(-\eta)=G(\eta)\circ\Sop$ gives $\Phi_{G(-\eta)}^1=\Sop\circ\Phi_{G(\eta)}^1\circ\Sop$.
Composing all factors proves (i).

(ii) We first prove the second identity $H^{(n-1)}(-\eta)=H^{(n-1)}(\eta)\circ\Sop$.
Because $I_j\circ\Sop=I_j$, part (i) gives $J_j^{(n-1)}(-\eta)=J_j^{(n-1)}(\eta)\circ\Sop$.
The coefficients $d_{ab}$ are even by \eqref{eq:equivariant-class} and $W\circ\Sop=-W$, so $(-\eta)\Kc^{(n-1)}(-\eta)=\eta\,\Kc^{(n-1)}(\eta)\circ\Sop$.
Finally, $H\circ\Sop=H$, and part (i) gives $\Psi^{(n-1)}(-\eta)=\Sop\circ\Psi^{(n-1)}(\eta)\circ\Sop$.
Hence
\[
\begin{aligned}
 H^{(n-1)}(-\eta)
 &=\big(H+(-\eta)\Kc^{(n-1)}(-\eta)\big)\circ\Psi^{(n-1)}(-\eta) \\
 &=\big(H+\eta\Kc^{(n-1)}(\eta)\big)\circ\Sop\circ\Sop\circ\Psi^{(n-1)}(\eta)\circ\Sop \\
 &=H^{(n-1)}(\eta)\circ\Sop .
\end{aligned}
\]
The first identity $H^{(n-1)}\circ\Rop=H^{(n-1)}$ follows in the same way from (i), $I_j\circ\Rop=I_j$, $W\circ\Rop=W$, and $H\circ\Rop=H$.
Both symmetries are linear and degree preserving, so they pass to the degree-$n$ component $\hat H_n$.

(iii) In the complex coordinates \eqref{eq:complex-convention}, $\Sop$ acts by $(q_1,p_1,q_2,p_2)\mapsto(q_1,p_1,-q_2,-p_2)$.
Since the second components satisfy $k_{q,2}+k_{p,2}\equiv k_{q,2}-k_{p,2}=\nu_2\pmod 2$, each monomial satisfies
\[
 \big(q^{k_q}p^{k_p}\big)\circ\Sop=(-1)^{\nu_2}\,q^{k_q}p^{k_p}.
\]
Comparing coefficients in $\hat H_n(-\eta)=\hat H_n(\eta)\circ\Sop$ gives $\hat h_{k_q,k_p}(-\eta)=(-1)^{\nu_2}\hat h_{k_q,k_p}(\eta)$, so $\hat h_{k_q,k_p}$ is odd whenever $\nu_2$ is odd.
Now let $n$ be even and $2(a+b)+2=n$, so that $I_1^aI_2^b\sqrt{I_1I_2}\cos\psi$ has degree $n$.
Since $\Sop$ fixes $I_1,I_2$ and sends $\psi$ to $\psi-\pi$, it changes the sign of this function.
Comparing its coefficients in $\hat H_n(-\eta)=\hat H_n(\eta)\circ\Sop$ gives $h_{ab}(-\eta)=-h_{ab}(\eta)$.
An odd formal series, or an odd real analytic function on $|\eta|\le\eta_0$, has the form \eqref{eq:hfactor}.

(iv) Suppose that $n$ is even.
By (ii) the known part $\hat H_n$ is $\Rop$-even, so Lemma~\ref{lem:solvability}(ii) leaves only the cosine direction in the first slow space associated with each $I_1^aI_2^b$.
By Lemma~\ref{lem:solvability}(i), the absorption equation splits into one equation for each pair $(a,b)$, and Lemma~\ref{lem:solvability}(iii) reduces each of them to
\[
h_{ab}+2\eta\sigma d_{ab}=0.
\]
At the formal level, $\sigma$ is even and $\sigma(0)=1$, so it is a unit of $\mathbb R[[\eta^2]]$.
At finite order, $\sigma$ is analytic and has a positive lower bound on every compact interval.
By (iii), $h_{ab}=\eta\widetilde h_{ab}(\eta^2)$, and cancelling the factor $\eta$ gives the unique solution \eqref{eq:dabsolution}, which is even in $\eta$.
Every $\Sm{\pm1}$ monomial has odd degree in each mode and hence even total degree, so for odd $n$ the absorption equation reduces to $0=0$ and no new counterterm coefficient occurs.

(v) The right-hand side of \eqref{eq:homological} contains only fast and higher slow monomials.
Dividing each coefficient by the eigenvalue \eqref{eq:eigenvalue} gives
\[
G_n=\sum\frac{-\hat h_{k_q,k_p}}{\langle k_p-k_q,i\Omega\rangle}\,q^{k_q}p^{k_p},
\]
where the sum is taken over the fast and higher slow monomials of degree $n$.
The fast divisors are invertible by Hypothesis~\ref{hyp:NR}, and the higher slow divisors are invertible by Corollary~\ref{prop:freq} and $\delta>0$.
The normalization condition sets all action and first slow components to zero, giving uniqueness.
Since the equation and the normalization condition are real, uniqueness also shows that $G_n$ is real.

The reflection $\Sop$ preserves each $\Sm{m}$ and $\Kmod$, while $\Rop$ maps $\Sm{m}$ to $\Sm{-m}$ and preserves $\Kmod$, so both commute with the projections $\Piker$ and $\Pislow$ and preserve the normalization condition \eqref{eq:gauge}.
For $\Sop$, since $\Omega$ is even in $\eta$ and $I_j\circ\Sop=I_j$, the operator $\mathcal L_0$ is the same at $\eta$ and $-\eta$ and satisfies $\mathcal L_0(G\circ\Sop)=(\mathcal L_0G)\circ\Sop$.
By (ii), $\hat H_n(-\eta)=\hat H_n(\eta)\circ\Sop$, so applying $\circ\,\Sop$ to \eqref{eq:homological} at parameter $\eta$ shows that $G_n(\eta)\circ\Sop$ solves the degree-$n$ homological equation at parameter $-\eta$ under the same normalization condition.
Uniqueness gives $G_n(-\eta)=G_n(\eta)\circ\Sop$.
For $\Rop$, since $\Rop$ is anti-symplectic and $(\Omega\cdot I)\circ\Rop=\Omega\cdot I$, we have $\mathcal L_0(G\circ\Rop)=-(\mathcal L_0G)\circ\Rop$.
Hence $-G_n\circ\Rop$ solves the same equation as $G_n$, with the same $\Rop$-even right-hand side by (ii), and uniqueness gives $G_n\circ\Rop=-G_n$.
\end{proof}

\subsection{Formal coupling normal form}\label{ssec:mainthm}

\begin{theorem}[coupling normal form at the formal level]
\label{thm:main}
Assume Hypothesis~\ref{hyp:abstract-symmetry} and the formal nonresonance conditions of Hypothesis~\ref{hyp:NR}.
For the coupling direction $W$, a fixed coupling scale $\kappa$, and the normal form style of Definition~\ref{def:style}, there are unique even counterterm coefficients $d_{ab}\in\mathbb R[[\eta^2]]$ and a unique formal symplectic map of the ordered form
\begin{equation}\label{eq:Psi-formal}
 \Psi
 =R_{\theta(\eta)}\circ\Phi_{G_3}^1\circ\Phi_{G_4}^1\circ\cdots,
\end{equation}
whose generators satisfy the normalization condition \eqref{eq:gauge}, such that
\begin{equation}\label{eq:NF}
 (H+\eta\Kc)\circ\Psi
 =H_{\rm nf}(I_1,I_2;\eta)
 =\Omega\cdot I
  +\sum_{a+b\ge2}A_{ab}(\eta)I_1^aI_2^b.
\end{equation}
Here $H_{\rm nf}$ denotes the resulting action normal form and $A_{ab}(\eta)$, $a,b\in\mathbb N_0$ with $a+b\ge2$, denote its action coefficients.
Moreover, $H_{\rm nf}$ is even in $\eta$, and
\begin{equation}\label{eq:dressedfactor}
 \Kc\circ\Psi
 =\Delta(I_1,I_2;\eta)(W\circ\Psi).
\end{equation}
\end{theorem}

\begin{proof}
For $n\ge2$ let $\mathrm A_n$ be the statement that there exist unique generators $G_3,\ldots,G_n$ satisfying the normalization condition \eqref{eq:gauge} and unique coefficients $d_{ab}$ with $2(a+b)+2\le n$ such that $\mathcal T_nH^{(n)}$ lies in the action kernel $\Kmod$.
The statement $\mathrm A_n$ also asserts that these generators and coefficients satisfy the equivariance conditions \eqref{eq:equivariant-class}.

The base case $\mathrm A_2$ holds by Proposition~\ref{thm:linear}, since $\mathcal T_2H^{(2)}=\Omega\cdot I$ and no generator or coefficient $d_{ab}$ with $a+b\ge1$ occurs.
Assume that $\mathrm A_{n-1}$ holds, and split $\mathcal T_nH^{(n)}\in\Kmod$ into $\mathcal T_{n-1}H^{(n)}\in\Kmod$ and the condition that the degree-$n$ component of $H^{(n)}$ lies in $\Kmod$.
By Proposition~\ref{prop:triangular}(ii), $\mathcal T_{n-1}H^{(n)}=\mathcal T_{n-1}H^{(n-1)}$, which does not involve $G_n$ or the new coefficients $d_{ab}$ with $2(a+b)+2=n$.
The first condition is therefore the condition in $\mathrm A_{n-1}$, and its unique solution fixes $G_3,\ldots,G_{n-1}$ and the coefficients with $2(a+b)+2\le n-1$.
With these fixed, the known part $\hat H_n$ is determined, and by Section~\ref{ssec:homological} the second condition becomes the homological and absorption equations \eqref{eq:homological}--\eqref{eq:absorption} for $G_n$ and the new coefficients.
The solution given by $\mathrm A_{n-1}$ satisfies the hypotheses of Proposition~\ref{prop:equivariant}.
For even $n$, part (iv) of that proposition gives the unique even new coefficients, and for odd $n$ no new coefficient occurs.
Part (v) gives the unique generator $G_n$ satisfying \eqref{eq:gauge}, and the equivariance conditions \eqref{eq:equivariant-class} then hold through degree $n$.
This proves $\mathrm A_n$.

Since $\Phi_{G_n}^1-\mathrm{id}$ starts at degree $n-1$, the Taylor jets of $\Psi^{(n)}$ stabilize as $n\to\infty$.
By Proposition~\ref{prop:triangular}(i), the jets of $J_j^{(n)}$ stabilize as well, and therefore the jets of $\Kc^{(n)}$ also stabilize.
These stable jets define the formal limits $\Psi$, $J_j$, and $\Kc$.
By Proposition~\ref{prop:triangular}(ii), $\mathcal T_nH^{(k)}=\mathcal T_nH^{(n)}$ for all $k\ge n$, so $\mathcal T_n\big((H+\eta\Kc)\circ\Psi\big)=\mathcal T_nH^{(n)}$.
By $\mathrm A_n$, this jet lies in $\Kmod$ and has quadratic part $\Omega\cdot I$.
Proposition~\ref{prop:equivariant}(ii) at degree $n+1$ then gives $H^{(n)}(-\eta)=H^{(n)}(\eta)\circ\Sop$.
This jet depends only on the actions and $I_j\circ\Sop=I_j$, so it is even in $\eta$.
Since $n$ is arbitrary, this gives \eqref{eq:NF} and the evenness of $H_{\rm nf}$, and $\mathrm A_n$ for all $n$ gives the uniqueness of $\Psi$ and of the even coefficients $d_{ab}$.
Finally, substituting \eqref{eq:dressedactions} into \eqref{eq:Heta} gives the factorization \eqref{eq:dressedfactor}.
\end{proof}

The factorization \eqref{eq:dressedfactor} gives the scalar \emph{recovery constraint}
\begin{equation}\label{eq:recovery}
 \eta\,\Delta(I_1,I_2;\eta)=0.
\end{equation}
On the tori satisfying \eqref{eq:recovery}, the auxiliary Hamiltonian term $\eta\Kc\circ\Psi$ vanishes, but its Hamiltonian vector field is in general nonzero.
The factor $\eta=0$ gives the original Hamiltonian $H$ and is called the \emph{trivial branch}.
At $\eta=0$ the rotation $R_{\theta(0)}$ is the identity.
By Proposition~\ref{prop:equivariant}(ii) and~(iii), the first slow coefficients of every known part also vanish there.
The recursion therefore reduces to the Birkhoff normalization of $H$ under \eqref{eq:gauge}, which yields the Birkhoff normal form $H_{\rm nf}(I;0)$ of $H$.
For $\eta\ne0$, the constraint reduces to the \emph{bifurcation equation} $\Delta(I_1,I_2;\eta)=0$, whose solutions form the \emph{nontrivial branch}.
This equation is the Hamiltonian counterpart of the equation $\check\Delta=0$ of \cite{LinChiba}.

\begin{remark}[content of the theorem and scope of uniqueness]
	\label{rem:uniqueness-scope}
	The factorization \eqref{eq:dressedfactor} is immediate from \eqref{eq:dressedactions}.
	The main content of the theorem is that the coupled construction exists and is unique degree by degree, by the triangularity of Proposition~\ref{prop:triangular} and the equivariant one-step solvability of Proposition~\ref{prop:equivariant}.
	Unlike a resonant normal form, \eqref{eq:NF} contains only action terms, and the first slow obstruction is carried by the counterterm coefficients $d_{ab}$.
	Uniqueness is relative to the coupling direction $W$, the coupling scale $\kappa$, the normal form style of Definition~\ref{def:style}, the normalization condition \eqref{eq:gauge}, and the ordered composition \eqref{eq:Psi-formal} of time-one flows, with one homogeneous generator per degree.
\end{remark}
% =====================================================================
\section{Finite-order dynamics and the bifurcation equation}
\label{sec:reduced}
% =====================================================================

This section truncates the construction of Section~\ref{sec:style} at a finite phase-space degree, where the normalizing map is analytic on a fixed neighborhood.
We estimate the remainder, construct approximately invariant tori and circles with their internal dynamics, and determine where the bifurcating branches leave the trivial branch.

\subsection{Finite-order normal form and remainder}\label{ssec:truncation}

For $d\ge1$ and $r>0$, write
\[
\mathcal D_r^{2d} :=\{z\in\mathbb C^{2d}:\|z\|_\infty<r\}, \qquad B_r^{2d}:=\mathcal D_r^{2d}\cap\mathbb R^{2d}.
\]
For a Hamiltonian $F$ analytic on $\mathcal D_{r_0}^{2d}$ and $0<r\le r_0$, set $\|X_F\|_r:=\sup_{z\in\mathcal D_r^{2d}}\|X_F(z)\|_\infty$.
We use $\mathbb T^2:=(\mathbb R/2\pi\mathbb Z)^2$ for the real angle torus.
Let $H_{\le N}:=\mathcal T_NH$ denote the truncation of the original Hamiltonian $H$ through phase-space degree $N$.
The following lemma is a parameter-uniform version of standard analytic estimates for Lie transforms \cite{Gio,GDFGS}.

\begin{lemma}[finite Lie transforms on nested analytic domains]
\label{lem:finite-lie}
Let $U\subset\mathbb C$ be a complex neighborhood of the interval $[-\eta_0,\eta_0]$.
Suppose that $F$ and the finite family $G_j$, $1\le j\le m$ with $m\in\mathbb N$, depend holomorphically on $(z,\eta)\in\mathcal D_{r_0}^{2d}\times U$ and that every $G_j$ vanishes to order at least three at the origin, uniformly for $|\eta|\le\eta_0$.
After decreasing $r_0$ if necessary, there are radii
\[
r_0>r_1>\cdots>r_m=:r_*>0
\]
such that, for $|t|\le1$ and $|\eta|\le\eta_0$, the complex Hamiltonian flow $\Phi_{G_j}^t$ exists on $\mathcal D_{r_j}^{2d}$ and maps it into $\mathcal D_{r_{j-1}}^{2d}$.
The finite composition $\Phi:=\Phi_{G_1}^1\circ\cdots\circ\Phi_{G_m}^1$ and its inverse are therefore holomorphic on fixed neighborhoods independent of $\eta$.

Let $R_{N+1}:=F\circ\Phi-\mathcal T_N(F\circ\Phi)$.
Then a constant $C_N$ exists, independent of $|\eta|\le\eta_0$ and $0<r\le r_*/2$, such that
\begin{equation}\label{eq:finite-lie-remainder}
 \|X_{R_{N+1}}\|_r\le C_Nr^N.
\end{equation}
The same conclusion holds when the composition is preceded or followed by an invertible linear map depending analytically on $\eta$, after changing the outer domain and the constant.
\end{lemma}

\begin{proof}
For $0<s<r_0$, uniform third-order vanishing gives a constant $C^{\rm flow}>0$ such that $\|X_{G_j}\|_s\le C^{\rm flow}s^2$ for every $j$ on a sufficiently small common polydisc.
Choose the nested radii so that $C^{\rm flow}r_{j-1}^2<r_{j-1}-r_j$.
Then, for every generator and $|t|\le1$, the complex ODE started in $\mathcal D_{r_j}^{2d}$ remains inside $\mathcal D_{r_{j-1}}^{2d}$, and continuation and holomorphic dependence give the asserted flows.
Since this chain does not depend on the order of the generators and $|t|\le1$ includes $t=-1$, it also gives the inverse $\Phi^{-1}=\Phi_{G_m}^{-1}\circ\cdots\circ\Phi_{G_1}^{-1}$.
Since $\mathcal T_NR_{N+1}=0$, we have $\mathcal T_{N-1}X_{R_{N+1}}=0$.
Cauchy estimates on $\mathcal D_{r_*}^{2d}$, followed by the Taylor estimate on $\mathcal D_r^{2d}$, give \eqref{eq:finite-lie-remainder}.
Compactness of $[-\eta_0,\eta_0]$ makes the bound uniform.
Such a linear map only changes the outer domain and the constant.
\end{proof}

\begin{proposition}[finite-order coupling normal form and remainder]
\label{prop:remainder}
Assume Hypothesis~\ref{hyp:abstract-symmetry}, fix $N\ge4$, and let $|\eta|\le\eta_0$ be an interval satisfying Hypothesis~\ref{hyp:NR}.
The recursion of Propositions~\ref{prop:triangular} and~\ref{prop:equivariant}, applied through degree $N$ with real analytic parameter dependence, gives even real analytic coefficients $d_{ab}(\eta)$ for $1\le a+b\le\lfloor(N-2)/2\rfloor$ and the stage-$N$ objects $\Psi^{(N)}$, $J_j^{(N)}$, $\Delta^{(N)}$, and $\Kc^{(N)}$ of \eqref{eq:partial-objects}--\eqref{eq:partial-Delta}.
The map $\Psi^{(N)}=R_{\theta(\eta)}\circ\Phi_{G_3}^1\circ\cdots\circ\Phi_{G_N}^1$ is real analytic and symplectic, and there exists a radius $r_{{\rm nf},N}>0$ such that
\begin{equation}\label{eq:finiteconj}
 (H_{\le N}+\eta\Kc^{(N)})\circ\Psi^{(N)}
 =H_{\rm nf}^{(N)}(I;\eta)+R_{N+1}(I,\varphi;\eta),
\end{equation}
where the finite-order action normal form $H_{\rm nf}^{(N)}$ is even in $\eta$ and $\mathcal T_NR_{N+1}=0$.
The finite-order counterterm factorizes exactly:
\begin{equation}\label{eq:finitefactor}
 \Kc^{(N)}\circ\Psi^{(N)}
 =\Delta^{(N)}(I_1,I_2;\eta)
  (W\circ\Psi^{(N)}).
\end{equation}
Moreover, for $0<r\le r_{{\rm nf},N}/2$,
\begin{equation}\label{eq:remainderfield}
 \sup_{|\eta|\le\eta_0}
 \|X_{R_{N+1}}\|_r\le C_Nr^N .
\end{equation}
\end{proposition}

\begin{proof}
Only finitely many analytic generators and divisors occur at fixed $N$.
Compactness of the parameter interval and Hypothesis~\ref{hyp:NR} give uniform coefficient bounds.
Applying Propositions~\ref{prop:triangular} and~\ref{prop:equivariant} through degree $N$, followed by Lemma~\ref{lem:finite-lie}, gives \eqref{eq:finiteconj}.
By Proposition~\ref{prop:equivariant}(v), $G_N$ also satisfies \eqref{eq:equivariant-class}, so $H^{(N)}(-\eta)=H^{(N)}(\eta)\circ\Sop$ follows as in the proof of Proposition~\ref{prop:equivariant}(ii).
Since $H_{\rm nf}^{(N)}=\mathcal T_NH^{(N)}$ by \eqref{eq:finiteconj} and $I_j\circ\Sop=I_j$, the normal form $H_{\rm nf}^{(N)}$ is even in $\eta$.
The identity $J_j^{(N)}\circ\Psi^{(N)}=I_j$ gives \eqref{eq:finitefactor}, and \eqref{eq:finite-lie-remainder} gives \eqref{eq:remainderfield}.
\end{proof}

\subsection{Invariant tori and their internal dynamics}\label{ssec:invariance}

On the tori satisfying $\Delta^{(N)}=0$, the auxiliary Hamiltonian term vanishes, and they become invariant tori of a Hamiltonian that differs from $H_{\le N}$ by the remainder.

\begin{proposition}[invariance with internal dynamics and residual estimate]
\label{prop:invariance}
Assume the setting of Proposition~\ref{prop:remainder}.
Let $I^*=(I_1^*,I_2^*)$ have positive components and satisfy
\begin{equation}\label{eq:constraint}
 \Delta^{(N)}(I^*;\eta)=0
\end{equation}
for some $|\eta|\le\eta_0$.
Suppose the torus $\{I=I^*\}$ is contained in $B_r^4$, with $r\le r_{{\rm nf},N}/2$, and define
\[
K_N(\varphi):=\Psi^{(N)}(I^*,\varphi), \qquad f_N(I,\varphi;\eta):=(W\circ\Psi^{(N)})(I,\varphi).
\]
The internal vector field on the angle torus $\mathbb T^2$ is
\begin{equation}\label{eq:innerflow}
 v_N(\varphi;I^*,\eta)
 :=\partial_I H_{\rm nf}^{(N)}(I^*;\eta)
 -\eta\,\partial_I\Delta^{(N)}(I^*;\eta)f_N(I^*,\varphi;\eta).
\end{equation}
When $I^*$ and $\eta$ are fixed, we abbreviate this field as $v_N(\varphi)$.
Then the recovered Hamiltonian
\[
H_{\rm rec}^{(N)} :=H_{\rm nf}^{(N)}\circ(\Psi^{(N)})^{-1}-\eta\Kc^{(N)} =H_{\le N}-R_{N+1}\circ(\Psi^{(N)})^{-1}
\]
leaves $K_N(\mathbb T^2)$ exactly invariant as a set and
\begin{equation}\label{eq:reconstructed-invariance}
 X_{H_{\rm rec}^{(N)}}(K_N)=DK_Nv_N.
\end{equation}
For the truncated original Hamiltonian $H_{\le N}$, there is a constant $\widetilde C_N>0$, independent of $r$ on the stated domain, such that the residual satisfies
\begin{equation}\label{eq:invariance-defect}
 \|X_{H_{\le N}}(K_N)-DK_Nv_N\|_{C^0}
 \le \widetilde C_Nr^N .
\end{equation}
Moreover, \eqref{eq:constraint} also holds for $-\eta$, and the corresponding tori satisfy
\[
K_N(\varphi;-\eta)=\Sop K_N(\varphi_1,\varphi_2+\pi;\eta), \qquad v_N(\varphi;I^*,-\eta)=v_N(\varphi_1,\varphi_2+\pi;I^*,\eta),
\]
so they are exchanged by $\Sop$ and have the same internal dynamics up to the shift $\varphi_2\mapsto\varphi_2+\pi$.
\end{proposition}

\begin{proof}
In normal-form coordinates,
\[
H_{\rm rec}^{(N)}\circ\Psi^{(N)} =H_{\rm nf}^{(N)}-\eta\Delta^{(N)}(I)f_N(I,\varphi).
\]
On the torus $I=I^*$, the constraint \eqref{eq:constraint} makes the action equations vanish and reduces the angular equations to \eqref{eq:innerflow}.
Combining \eqref{eq:finiteconj} and \eqref{eq:finitefactor} shows that the residual for $H_{\le N}$ is the push-forward of $X_{R_{N+1}}$.
Uniform bounds on $D\Psi^{(N)}$ and \eqref{eq:remainderfield} prove \eqref{eq:invariance-defect}.

Since $\Delta^{(N)}$ is even in $\eta$, \eqref{eq:constraint} also holds for $-\eta$.
By Proposition~\ref{prop:equivariant}(v), $G_N$ also satisfies \eqref{eq:equivariant-class}, so $\Psi^{(N)}(-\eta)=\Sop\circ\Psi^{(N)}(\eta)\circ\Sop$ follows as in the proof of Proposition~\ref{prop:equivariant}(i).
Since $\Sop$ acts on the normal-form action--angle variables by $\varphi_2\mapsto\varphi_2+\pi$, this gives the identity for $K_N$.
For the internal field, $H_{\rm nf}^{(N)}$ and $\Delta^{(N)}$ are even in $\eta$ by Proposition~\ref{prop:remainder}, while $W\circ\Sop=-W$ gives $f_N(I^*,\varphi;-\eta)=-f_N(I^*,\varphi_1,\varphi_2+\pi;\eta)$, which proves the identity for $v_N$.
\end{proof}

\begin{remark}[tangential correction from the auxiliary Hamiltonian term]
\label{rem:tangential}
Let $e_j$ be the $j$-th coordinate vector on the angle torus.
Since $\Psi^{(N)}$ is symplectic and $J_j^{(N)}\circ\Psi^{(N)}=I_j$, the Hamiltonian vector fields of $J_j^{(N)}$ satisfy $X_{J_j^{(N)}}\circ K_N=DK_N\,e_j$.
The product rule and \eqref{eq:constraint} therefore give
\[
 X_{\eta\Kc^{(N)}}\circ K_N
 =DK_N\bigl[\eta f_N(I^*,\varphi;\eta)
       \,\partial_I\Delta^{(N)}(I^*;\eta)\bigr].
\]
Thus the auxiliary Hamiltonian term $\eta\Kc^{(N)}$ vanishes on $K_N(\mathbb T^2)$, while its vector field is tangent to $K_N(\mathbb T^2)$ and in general nonzero.
Subtracting this field in $H_{\rm rec}^{(N)}$ gives the negative correction in \eqref{eq:innerflow}.
Since $\partial_I\Delta^{(N)}(I^*;\eta)$ is in general not parallel to the frequency vector $\partial_IH_{\rm nf}^{(N)}(I^*;\eta)$, this tangential correction is not merely a reparametrization of time.
\end{remark}

\begin{remark}[angular straightening]
\label{prop:rotation}
The field $v_N$ is in general angle dependent, so \eqref{eq:invariance-defect} is not yet a constant-frequency invariance equation $X_H(K)-DK\boldsymbol\rho=0$.
If a diffeomorphism $h$ of $\mathbb T^2$ and a constant frequency vector $\boldsymbol\rho$ satisfy
\begin{equation}\label{eq:straightening}
 Dh(\vartheta)\boldsymbol\rho=v_N(h(\vartheta)),
\end{equation}
then $D(K_N\circ h)\boldsymbol\rho=(DK_N\circ h)\,Dh\,\boldsymbol\rho$ shows that $K_N\circ h$ satisfies \eqref{eq:invariance-defect} with the same constant, relative to $\boldsymbol\rho$.
Solving \eqref{eq:straightening} is a small-divisor problem on $\mathbb T^2$ and is not addressed here.
\end{remark}

Proposition~\ref{prop:invariance} assumes $I_1^*,I_2^*>0$.
On the axis $I_2=0$, however, the angle $\varphi_2$ is not defined.
The next proposition treats the axis directly in the Cartesian canonical variables $(Q,P)$, where no angle coordinate is needed for the vanishing mode.

\begin{proposition}[axis circles in canonical variables]
\label{prop:axis}
Assume the setting of Proposition~\ref{prop:remainder}.
Let $I_1^*>0$ satisfy
\begin{equation}\label{eq:axis-constraint}
 \Delta^{(N)}(I_1^*,0;\eta)=0
\end{equation}
for some $|\eta|\le\eta_0$, and let
\[
\mathcal O:=\big\{(Q_1,P_1,0,0):\tfrac12(Q_1^2+P_1^2)=I_1^*\big\} \subset B_r^4, \qquad r\le r_{{\rm nf},N}/2,
\]
parameterized by $\mathcal O(\varphi_1):=(\sqrt{2I_1^*}\cos\varphi_1,-\sqrt{2I_1^*}\sin\varphi_1,0,0)$ for $\varphi_1\in\mathbb T^1:=\mathbb R/2\pi\mathbb Z$.
Define
\[
K_N^{\rm ax}(\varphi_1):=\Psi^{(N)}(\mathcal O(\varphi_1)), \qquad f_N^{\rm ax}(\varphi_1;I_1^*,\eta) :=\big(W\circ\Psi^{(N)}\big)(\mathcal O(\varphi_1)).
\]
Then:
\begin{itemize}
\item[(i)] the circle $\mathcal O$ and its image $K_N^{\rm ax}(\mathbb T^1)$ are invariant under the flows of $H_{\rm rec}^{(N)}\circ\Psi^{(N)}$ and $H_{\rm rec}^{(N)}$, respectively.
\item[(ii)] the motion on $\mathcal O$ is $\dot\varphi_1=v_N^{\rm ax}(\varphi_1)$ with
\begin{equation}\label{eq:axis-field}
 v_N^{\rm ax}(\varphi_1;I_1^*,\eta)
 :=\partial_{I_1}H_{\rm nf}^{(N)}(I_1^*,0;\eta)
 -\eta\,\partial_{I_1}\Delta^{(N)}(I_1^*,0;\eta)\,
 f_N^{\rm ax}(\varphi_1;I_1^*,\eta).
\end{equation}
With $I_1^*$ and $\eta$ fixed, we abbreviate these axis functions as $f_N^{\rm ax}(\varphi_1)$ and $v_N^{\rm ax}(\varphi_1)$.
There is a radius $r_{\rm ax}\in(0,r_{{\rm nf},N}/2]$, uniform for $|\eta|\le\eta_0$, such that $v_N^{\rm ax}\ge\tfrac12\Omega_1(\eta)>0$ on $\mathbb T^1$ whenever $\mathcal O\subset B_{r_{\rm ax}}^4$.
The circle $\mathcal O$ is therefore a closed orbit of $H_{\rm rec}^{(N)}\circ\Psi^{(N)}$ with period $T:=\int_0^{2\pi}v_N^{\rm ax}(\varphi_1)^{-1}\,d\varphi_1$.
\item[(iii)] the residual estimate holds on the axis with the same constant $\widetilde C_N$ as in \eqref{eq:invariance-defect},
\begin{equation}\label{eq:axis-defect}
 \big\|X_{H_{\le N}}(K_N^{\rm ax})
 -DK_N^{\rm ax}\,v_N^{\rm ax}\big\|_{C^0}
 \le\widetilde C_Nr^N .
\end{equation}
\item[(iv)] \eqref{eq:axis-constraint} also holds for $-\eta$, and the corresponding circles satisfy
\[
K_N^{\rm ax}(\varphi_1;-\eta)=\Sop K_N^{\rm ax}(\varphi_1;\eta), \qquad v_N^{\rm ax}(\varphi_1;I_1^*,-\eta)=v_N^{\rm ax}(\varphi_1;I_1^*,\eta),
\]
so they are exchanged by $\Sop$ and have identical tangential speed.
\end{itemize}
\end{proposition}

\begin{proof}
(i) On the plane $\{Q_2=P_2=0\}$, the $(Q_2,P_2)$-gradient of any function of the actions vanishes, because $\nabla_{(Q_2,P_2)}I_2=(Q_2,P_2)$ and $\nabla_{(Q_2,P_2)}I_1=0$.
Hence, on this plane, the $(Q_2,P_2)$-components of the Hamiltonian vector field of $H_{\rm rec}^{(N)}\circ\Psi^{(N)}$ reduce to the single term $-\eta\,\Delta^{(N)}(I_1,0;\eta)\,\mathbb J_2\nabla_{(Q_2,P_2)} \big(W\circ\Psi^{(N)}\big)$, where $\mathbb J_2$ is the $2\times2$ symplectic matrix.
On $\mathcal O$ this term carries the factor \eqref{eq:axis-constraint} and vanishes.
For the radial direction, $\{I_1,H_{\rm nf}^{(N)}(I)\}=\{I_1,\Delta^{(N)}(I)\}=0$, so
\[
\dot I_1 =-\eta\,\Delta^{(N)}(I)\,\{I_1,W\circ\Psi^{(N)}\},
\]
which also vanishes on $\mathcal O$.
Hence $\mathcal O$ is invariant, and conjugation by $\Psi^{(N)}$ transports the invariance to $K_N^{\rm ax}(\mathbb T^1)$.

(ii) Since $I_1^*>0$, the angle $\varphi_1$ is a regular coordinate for the first mode near $\mathcal O$.
With the sign convention of \eqref{eq:diagonal-actions}, the tangential equation is therefore $\dot\varphi_1=\partial_{I_1}\big(H_{\rm rec}^{(N)}\circ\Psi^{(N)}\big)\big|_{\mathcal O}$, which is \eqref{eq:axis-field}.
By construction $\partial_{I_1}H_{\rm nf}^{(N)}(I_1^*,0;\eta)=\Omega_1(\eta)+O(r^2)$.
Moreover $W\circ\Psi^{(N)}$ has no constant or linear part, so $f_N^{\rm ax}=O(r^2)$, and $\partial_{I_1}\Delta^{(N)}(I_1^*,0;\eta)$ is bounded for $|\eta|\le\eta_0$.
Hence $v_N^{\rm ax}=\Omega_1(\eta)+O(r^2)$ uniformly for $|\eta|\le\eta_0$.
Since $\Omega_1(\eta)\ge\omega_1>0$ by Corollary~\ref{prop:freq}, a uniform $r_{\rm ax}$ with the asserted lower bound exists.
Then $v_N^{\rm ax}$ is positive on $\mathbb T^1$, so the solution of $\dot\varphi_1=v_N^{\rm ax}(\varphi_1)$ moves monotonically around $\mathcal O$.
It returns to its initial point after the time $T$.

(iii) By \eqref{eq:finiteconj}, $H_{\le N} =H_{\rm rec}^{(N)}+R_{N+1}\circ(\Psi^{(N)})^{-1}$.
Part (i) and \eqref{eq:axis-field} give $X_{H_{\rm rec}^{(N)}}(K_N^{\rm ax}) =DK_N^{\rm ax}v_N^{\rm ax}$, so the left-hand side of \eqref{eq:axis-defect} is the push-forward of $X_{R_{N+1}}$ by the symplectic map $\Psi^{(N)}$.
The bound \eqref{eq:remainderfield} and the uniform bound on $D\Psi^{(N)}$ give the estimate.

(iv) Since $\Delta^{(N)}$ is even in $\eta$, \eqref{eq:axis-constraint} also holds for $-\eta$, and as in the proof of Proposition~\ref{prop:invariance}, $\Psi^{(N)}(-\eta)=\Sop\circ\Psi^{(N)}(\eta)\circ\Sop$.
Since $\Sop\,\mathcal O(\varphi_1)=\mathcal O(\varphi_1)$, this gives $K_N^{\rm ax}(\varphi_1;-\eta)=\Sop K_N^{\rm ax}(\varphi_1;\eta)$.
For the tangential field, $H_{\rm nf}^{(N)}$ and $\Delta^{(N)}$ are even in $\eta$ by Proposition~\ref{prop:remainder}, while $W\circ\Sop=-W$ gives $f_N^{\rm ax}(\varphi_1;I_1^*,-\eta)=-f_N^{\rm ax}(\varphi_1;I_1^*,\eta)$, which proves the identity for $v_N^{\rm ax}$.
\end{proof}

\begin{remark}[no straightening obstruction on the axis]
\label{rem:axis-straightening}
On the axis the straightening of Remark~\ref{prop:rotation} always exists, because the invariant circle is one-dimensional and $v_N^{\rm ax}$ is positive by Proposition~\ref{prop:axis}(ii).
Set $\rho_{\rm ax}:=2\pi/T$.
The circle diffeomorphism $h:\mathbb T^1\to\mathbb T^1$ defined implicitly by
\[
\vartheta=\rho_{\rm ax}\int_0^{h(\vartheta)} \frac{ds}{v_N^{\rm ax}(s)}
\]
satisfies $h'(\vartheta)\rho_{\rm ax}=v_N^{\rm ax}(h(\vartheta))$, and $h(\vartheta+2\pi)=h(\vartheta)+2\pi$.
Hence $K_N^{\rm ax}\circ h$ satisfies the axis residual estimate \eqref{eq:axis-defect} relative to the constant frequency $\rho_{\rm ax}$.
The straightening problem of Remark~\ref{prop:rotation} therefore remains only for the approximately invariant tori.
\end{remark}

\subsection{Finite-order bifurcation equation}\label{ssec:generic-fourth}

This subsection determines, at finite order, where the nontrivial branch meets the trivial branch.
For $N\ge4$, Proposition~\ref{prop:remainder} shows that the coefficients of $\Delta^{(N)}(I;\eta)$ are even and real analytic on every admissible interval.
Therefore there is a unique function $F_N(I;u)$, polynomial in $I$ and analytic in $u$ near $u=0$, such that
\begin{equation}\label{eq:F-u}
 F_N(I;\eta^2)=\Delta^{(N)}(I;\eta),
 \qquad
 F_N(I;u)=\ell_N(I)+u\,\zeta_N(I)+O(u^2),
\end{equation}
where
\begin{equation}\label{eq:LQ}
 \ell_N(I):=\Delta^{(N)}(I;0),
 \qquad
 \zeta_N(I):=\partial_uF_N(I;0).
\end{equation}

\begin{theorem}[critical curve and bifurcating branches]
\label{thm:line}
Fix $N\ge4$, let $I^0$ lie on the \emph{critical curve} $\ell_N(I)=0$, and assume the nondegeneracy condition $\zeta_N(I^0)\ne0$.
Then there is a unique local analytic root branch $u=u(I)$ of $F_N(I;u)=0$ with
\begin{equation}\label{eq:localu}
 u(I)=-\frac{\ell_N(I)}{\zeta_N(I^0)}+O(|I-I^0|^2).
\end{equation}
For every $I$ with $u(I)>0$, the values $\eta=\pm\sqrt{u(I)}$ are symmetry-related real roots of the order-$N$ bifurcation equation $\Delta^{(N)}(I;\eta)=0$.
They form the two \emph{bifurcating branches}.
\end{theorem}

\begin{proof}
Since $F_N(I^0;0)=\ell_N(I^0)=0$ and $\partial_uF_N(I^0;0)=\zeta_N(I^0)\ne0$, the analytic implicit function theorem gives a unique local root branch $u(I)$ with $u(I^0)=0$.
Expanding $F_N(I;u(I))=0$ at $I^0$ gives \eqref{eq:localu}.
The last statement is the two-to-one real change of parameter $u=\eta^2$.
\end{proof}

At the lowest order,
\begin{equation}\label{eq:Delta4}
 \Delta^{(4)}(I;\eta)=d_{00}+d_{10}(\eta)I_1+d_{01}(\eta)I_2,
 \qquad
 \ell_4(I)=d_{00}+d_{10}(0)I_1+d_{01}(0)I_2,
\end{equation}
so the critical curve of order four is the \emph{critical line} $\ell_4(I)=0$.
The coefficients $d_{ab}$ do not change with $N$ by Proposition~\ref{prop:triangular}(ii), and the terms with $a+b\ge2$ are at least quadratic in the actions.
Hence $\ell_N=\ell_4+O(|I|^2)$, and the critical line approximates the critical curves of all orders at small actions.
If $d_{10}(0)>0$, the critical line meets the axis $I_2=0$ at the threshold
\begin{equation}\label{eq:threshold}
 I_1^{\min}=-\frac{d_{00}}{d_{10}(0)}=\frac{\kappa\delta}{d_{10}(0)},
\end{equation}
which is small when the frequency difference is small.

\begin{remark}[pitchfork bifurcation of the recovery constraint]
\label{rem:critical}
Because $\Delta^{(N)}$ is even in $\eta$, near $(I^0,0)$ the roots of $\eta\,\Delta^{(N)}(I;\eta)=0$ undergo a $\mathbb Z_2$-symmetric pitchfork bifurcation in $\eta$, with $\ell_N(I)$ playing the role of the bifurcation parameter.
However, this bifurcation concerns the finite-order recovery constraint, not the dynamics.
Along the bifurcating branches, Propositions~\ref{prop:invariance} and~\ref{prop:axis} provide approximately invariant tori and circles of the original system, and their family bifurcates from the trivial branch in this approximate sense.
A dynamical bifurcation would require in addition the persistence of these objects and a Floquet or variational degeneracy of the family on the trivial branch at the critical curve, both of which lie outside this paper.
\end{remark}

\begin{remark}[the coupling parameter on the bifurcating branches]
\label{rem:order-parameter}
For a fixed coupling scale $\kappa$, $\eta$ is a signed branch coordinate on the bifurcating branches.
Its value determines the principal rotation angle, its sign distinguishes the two branches related by $\Sop$, and its square $u=\eta^2$ is the reflection-invariant branch variable.
Since $R_{\theta(\eta)}$ is the linear part of the normalizing map, a nonzero $\eta$ mixes the two modes, so $\eta$ measures how far the linear part of the bifurcating tori and circles departs from the uncoupled modes of the trivial branch.
On the bifurcating branches, $\eta$ is determined by the actions, grows from zero across the critical curve, and is analogous to a Lyapunov--Schmidt amplitude.
The critical curves do not depend on $\kappa$ and are therefore intrinsic.
\end{remark}

% =====================================================================
\section{Application to the spatial CRTBP}\label{sec:crtbp}
% =====================================================================

This section applies the theory of Sections~\ref{sec:abstract}--\ref{sec:reduced} to the spatial CRTBP.
We reduce the problem to the center manifold, verify the abstract hypotheses, specialize the bifurcation equation, reconstruct the halo orbits and quasi-halo tori, and conclude with a numerical example for the Sun--Earth $L_1$.

\subsection{Dynamical model at collinear points}\label{ssec:setting}

Let $0<\mu\le\tfrac12$ be the mass parameter and let $L_i$ be the collinear point $L_1$ or $L_2$.
We use the standard scaled synodic coordinates $(x,y,z,p_x,p_y,p_z)$ centered at $L_i$.
The Hamiltonian is
\begin{equation}\label{eq:H}
 \mathcal H
 =\tfrac12(p_x^2+p_y^2+p_z^2)+y p_x-xp_y
 -\sum_{n\ge2}c_n(\mu)\rho^nP_n(x/\rho),
\end{equation}
where $\rho:=\sqrt{x^2+y^2+z^2}$, $P_n$ is the Legendre polynomial of degree $n$, and the coefficients $c_n(\mu)$ of the Legendre expansion of the potential about $L_i$ are given in closed form in \cite{JM}.
The anti-symplectic reversor $(x,y,z,p_x,p_y,p_z)\mapsto(x,-y,z,-p_x,p_y,-p_z)$ and the symplectic vertical reflection $(z,p_z)\mapsto(-z,-p_z)$ leave $\mathcal H$ invariant.

The linearization at $L_i$ has the hyperbolic eigenvalue $\lambda$ and the planar and vertical frequencies $\omega_1$, $\omega_2$ given by
\begin{equation}\label{eq:frequencies}
 \omega_1^2=\frac{2-c_2+\sqrt{9c_2^2-8c_2}}{2},
 \qquad
 \omega_2=\sqrt{c_2},
 \qquad
 \lambda^2=\frac{c_2-2+\sqrt{9c_2^2-8c_2}}{2}.
\end{equation}
At each collinear point and for every mass parameter $0<\mu\le\tfrac12$ one has $c_2>1$, so the equilibrium is of saddle$\times$center$\times$center type \cite{JM}.
Squaring reduces $\omega_1^2>c_2$ to $4c_2>4$, so $\omega_1>\omega_2=\sqrt{c_2}>1$.
The frequency difference and the frequency sum of \eqref{eq:abstract-frequencies} are
\begin{equation}\label{eq:detuning}
 \delta=\omega_1-\omega_2>0,
 \qquad
 g=\omega_1+\omega_2 .
\end{equation}

For every mass parameter the frequency difference is small compared with the frequency sum.
The collinear points thus lie in the near-resonance regime.
The bifurcation threshold \eqref{eq:threshold} is proportional to $\delta$, so the bifurcating branches meet the axis $I_2=0$ at small amplitude.

The quadratic part of $\mathcal H$ is
\[
\mathcal H_2=\tfrac12(p_x^2+p_y^2)+yp_x-xp_y-c_2x^2 +\tfrac{c_2}{2}y^2+\tfrac12p_z^2+\tfrac{c_2}{2}z^2.
\]

\begin{lemma}[symplectic normalization of $\mathcal H_2$]\label{lem:linlayer}
There is a real linear symplectic map
\[
L:(\bQ_h,\bP_h,\bQ_1,\bP_1,\bQ_2,\bP_2) \longmapsto(x,y,z,p_x,p_y,p_z)
\]
such that
\begin{equation}\label{eq:H2diag}
 \mathcal H_2\circ L
 =\lambda \bQ_h\bP_h
 +\frac{\omega_1}{2}(\bQ_1^2+\bP_1^2)
 +\frac{\omega_2}{2}(\bQ_2^2+\bP_2^2).
\end{equation}
The pairs $(\bQ_1,\bP_1)$ and $(\bQ_2,\bP_2)$ are the planar and vertical elliptic pairs, respectively.
The normalization constants are defined successively by
\begin{equation}\label{eq:mxmz}
 \xi=-\frac{\omega_1^2+1+2c_2}{2\omega_1}<0,
 \qquad
 s_\omega:=\omega_1(1+\xi^2)+2\xi,
 \qquad
 \mx:=\sqrt{\frac{2}{s_\omega}},
 \qquad
 \mz:=\sqrt{\frac{2}{\omega_2}}.
\end{equation}
In the reference action--angle variables \eqref{eq:reference-actions}, the physical coordinates $x$ and $z$ then satisfy at linear order
\begin{equation}\label{eq:xzrows0}
 x=\mx\sqrt{\bar I_1}\,\cos\bar\varphi_1+(\text{hyperbolic}),
 \qquad
 z=\mz\sqrt{\bar I_2}\,\cos\bar\varphi_2.
\end{equation}
In these coordinates the reversor and the vertical reflection preserve the center subspace $\{\bQ_h=\bP_h=0\}$.
Their restrictions to it act in the coordinates $(\bQ_1,\bP_1,\bQ_2,\bP_2)$ as the reversor $\Rop$ of \eqref{eq:abstract-R} and the reflection $\Sop$ of \eqref{eq:abstract-S}.
\end{lemma}

The construction of $L$ is given in Appendix~\ref{app:linnorm}.

\subsection{Reduction to the center manifold}
\label{ssec:cm}

Let $\iota:\mathbb R^4\to\mathbb R^6$, $\iota(\bQ_1,\bP_1,\bQ_2,\bP_2):=(0,0,\bQ_1,\bP_1,\bQ_2,\bP_2)$, be the map that places the center variables on the center subspace $\mathcal C:=\{\bQ_h=\bP_h=0\}$.
The construction order is
\begin{equation}\label{eq:architecture}
 \mathcal H
 \xrightarrow[\text{from }\mathbb R^6\text{ to }\mathbb R^4]{\text{center manifold reduction}}
 H_{\le N}
 \xrightarrow[\text{on }\mathbb R^4]{\text{auxiliary family}}
 H_{\le N}+\eta\Kc^{(N)}.
\end{equation}

\begin{lemma}[finite-order center manifold reduction and reconstruction]
\label{lem:cm}
For every $N\ge3$ there are polynomial generators $\Gamma_3,\ldots,\Gamma_N$ such that, with the local symplectic map
\[
\Phi_{\rm cm}^{(N)} :=\Phi_{\Gamma_3}^1\circ\cdots\circ\Phi_{\Gamma_N}^1,
\]
the degree-$N$ truncation
\begin{equation}\label{eq:HcmN}
 \mathcal H^{\rm cm}_{\le N}
 :=\mathcal T_N\big((\mathcal H\circ L)
   \circ\Phi_{\rm cm}^{(N)}\big)
\end{equation}
depends on the hyperbolic variables only through the product $\bQ_h\bP_h$.
Hence $\mathcal C$ is an exactly invariant center manifold of $\mathcal H^{\rm cm}_{\le N}$ and
\begin{equation}\label{eq:HN}
 H_{\le N}:=\mathcal H^{\rm cm}_{\le N}\circ\iota
\end{equation}
is a two-degree-of-freedom polynomial Hamiltonian with quadratic part $H_2=\tfrac{\omega_1}{2}(\bQ_1^2+\bP_1^2)+\tfrac{\omega_2}{2}(\bQ_2^2+\bP_2^2)$.

For a suitable radius $r_{{\rm cm},N}>0$, define
\begin{equation}\label{eq:EN}
 E_N:=L\circ\Phi_{\rm cm}^{(N)}\circ\iota
 :B_{r_{{\rm cm},N}}^4\longrightarrow\mathbb R^6
\end{equation}
and the center-manifold residual
\begin{equation}\label{eq:cmdefect}
 e_{{\rm cm},N}
 :=X_{\mathcal H}\circ E_N-DE_NX_{H_{\le N}}.
\end{equation}
Then there is a constant $C_{{\rm cm},N}>0$ such that, for $0<r\le r_{{\rm cm},N}/2$,
\begin{equation}\label{eq:cmdefectbound}
 \|e_{{\rm cm},N}\|_{C^0(B_r^4)}
 \le C_{{\rm cm},N}r^N.
\end{equation}
If each generator has no components on monomials with equal hyperbolic exponents, the generators are unique and $\Gamma_n$ does not depend on $N$.
Then $H_{\le N}$ is invariant under $\Rop$ and $\Sop$, $E_N$ intertwines them with the reversor and the vertical reflection of $\mathcal H$, and $H_{\le N}=\mathcal T_NH$ for a formal reduced Hamiltonian $H$, as in Section~\ref{ssec:truncation}.
\end{lemma}

The proof is given in Appendix~\ref{app:cm}.

\subsection{Verification of the abstract hypotheses}
\label{ssec:verification}

This subsection verifies \eqref{eq:abstract-frequencies} and Hypothesis~\ref{hyp:abstract-symmetry} for the reduced Hamiltonian, identifies the coupling direction of Definition~\ref{def:abstract-W} with the first slow projection of the product $xz$, and specializes the coupling scale \eqref{eq:abstract-d00}.

The product $xz$ of the physical coordinates is odd under the vertical reflection and therefore does not occur in the CRTBP Hamiltonian \eqref{eq:H}.
It is the physical coupling function that produces the coupling term of the modified vertical equation in \cite{LinChiba}.

By \eqref{eq:xzrows0}, the quadratic part of $xz\circ L$ restricted to $\mathcal C$ is
\begin{equation}\label{eq:xzdecomp}
 (xz\circ L)_2|_{\mathcal C}
 =\chixz\sqrt{\bar I_1\bar I_2}\,\cos(\bar\varphi_1-\bar\varphi_2)
 +\chixz\sqrt{\bar I_1\bar I_2}\,\cos(\bar\varphi_1+\bar\varphi_2),
\end{equation}
with the coefficient
\begin{equation}\label{eq:chixzdef}
 \chixz:=\tfrac12\mx\mz
 =\frac{1}{\sqrt{s_\omega\omega_2}}>0.
\end{equation}
The first term is first slow, with divisor $\delta$, and equals $\tfrac12\chixz W$ for the coupling direction $W$ of \eqref{eq:abstract-W}.
The second term is fast, with divisor $g$.
Thus $W$ is the normalized first slow projection of $(xz\circ L)_2|_{\mathcal C}$.
Since $\chixz>0$, the product $xz$ projects nontrivially onto $W$.

The auxiliary family uses $W$, not the full product $xz$, whose fast and hyperbolic components are discarded.
Adding the full product to $\mathcal H$ would change the linear frequencies and the center manifold reduction, and would therefore define a different scheme.

To match \cite{LinChiba}, we choose the coupling scale
\begin{equation}\label{eq:dictionary}
 \kappa=g,
 \qquad\text{so that}\qquad
 d_{00}=-g\delta=\omega_2^2-\omega_1^2 .
\end{equation}

\begin{corollary}[CRTBP realization of the abstract hypotheses]
\label{cor:crtbp-verification}
In the reference coordinates $(\bQ_1,\bP_1,\bQ_2,\bP_2)$ of $\mathcal C$, the reduced Hamiltonian $H$ satisfies \eqref{eq:abstract-frequencies} and Hypothesis~\ref{hyp:abstract-symmetry}, with the frequencies \eqref{eq:frequencies} and the reversor and reflection of Lemma~\ref{lem:linlayer}.
With the coupling direction $W$ of Definition~\ref{def:abstract-W} and the coupling scale \eqref{eq:dictionary}, the formal results of Sections~\ref{sec:linear}--\ref{sec:reduced} hold under the formal nonresonance conditions of Hypothesis~\ref{hyp:NR}, and the finite-order results hold on admissible intervals.
\end{corollary}

\begin{proof}
Lemma~\ref{lem:cm} gives the reduced Hamiltonian and its inherited symmetries.
Lemma~\ref{lem:linlayer} identifies the two elliptic modes and the action of $\Rop$ and $\Sop$, and the frequencies satisfy \eqref{eq:abstract-frequencies} by \eqref{eq:detuning}.
\end{proof}

The admissible intervals in Corollary~\ref{cor:crtbp-verification} are explicit.
Since $\langle\nu,\Omega\rangle=\tfrac12(\nu_1+\nu_2)g+\tfrac12(\nu_1-\nu_2)\delta/\sigma(\eta)$, and the angular index of a fast monomial of degree at most $N$ has $\nu_1+\nu_2\ne0$ and $|\nu_1-\nu_2|\le N$, the interval $|\eta|\le\eta_0$ is admissible whenever $N\delta/\sigma(\eta_0)<g$.
The formal nonresonance conditions hold exactly when $\omega_1/\omega_2$ is irrational, which is the case for all but countably many mass parameters.

\begin{remark}[relation to the coupling coefficient of the modified equation]
\label{rem:LCdict}
The modified linear vertical equation $\ddot z+c_2z=\check\eta\,d_{00}x$ of \cite{LinChiba} is the vertical equation of the quadratic Hamiltonian $\mathcal H_2-\check\eta\,d_{00}\,xz$, with the same scale $d_{00}=\omega_2^2-\omega_1^2$.
The same term also adds $\check\eta\,d_{00}z$ to the planar equation for $x$.
By \eqref{eq:xzdecomp}, its first slow component on the center manifold is $-\tfrac12\chixz\,\check\eta\,d_{00}W$.
The auxiliary quadratic family $H_2(\eta)$ retains only this component and discards the fast one, so the first slow components of the two quadratic Hamiltonians agree under
\[
\eta=-\tfrac12\chixz\,\check\eta.
\]
In \cite{LinChiba,LLC}, $\check\eta$ is the coefficient of the induced planar harmonic in
\[
z=\beta\cos\check\theta_2+\check\eta\,\alpha\cos\check\theta_1+\cdots
\]
and is read as an out-of-plane to in-plane amplitude ratio, whereas the present $\eta$ is the parameter of the auxiliary family that produces such a ratio.
Applying $R_{\theta(\eta)}$ to \eqref{eq:xzrows0}, with $\kappa=g$ in \eqref{eq:rotation-angle}, shows that the eigenmode continuing the planar reference mode has the linearized amplitude ratio
\[
\frac{\mz}{\mx}\,\big|\tan\theta(\eta)\big|
=\frac{\mz}{\mx}\,g\,|\eta|+O(|\eta|^{3})
=\frac{g}{2\omega_2}\,|\check\eta|+O(|\check\eta|^3),
\]
where the last equality uses the correspondence above.
The solution of \cite{LinChiba} has the ratio $|\check\eta|$ instead.
The factor $g/(2\omega_2)$ arises because discarding the fast component replaces the forcing denominator $\omega_2^2-\omega_1^2$ by $2\omega_2(\omega_2-\omega_1)$.
Thus the scale of $\eta$ is set by $W$ and $\kappa$, not by a normalization of the solution (Remark~\ref{rem:eta-status}).
\end{remark}

\subsection{Bifurcation equation in the CRTBP}\label{ssec:d4}

By Corollary~\ref{cor:crtbp-verification}, Theorem~\ref{thm:line} applies with $I_1$ and $I_2$ the planar and vertical actions.

\begin{corollary}[critical curves in the CRTBP]\label{thm:line-crtbp}
Let $N\ge4$. At every point $I^0$ of the critical curve $\ell_N(I)=0$ with $\zeta_N(I^0)\ne0$, the order-$N$ bifurcation equation has the unique local root branch \eqref{eq:localu}, and for every $I$ with $u(I)>0$ it has the symmetry-related real roots $\eta=\pm\sqrt{u(I)}$.
For $N=4$ the critical curve is the critical line
\begin{equation}\label{eq:line}
 d_{10}(0)I_1+d_{01}(0)I_2=g\delta .
\end{equation}
If $d_{10}(0)>0$, it meets the axis $I_2=0$ at the bifurcation threshold \eqref{eq:threshold} with $\kappa=g$.
If $d_{01}(0)\le0$, it does not meet the vertical axis $I_1=0$, so no bifurcating branches emerge from that axis.
\end{corollary}

\begin{proof}
The first statement is Theorem~\ref{thm:line} in the setting of Corollary~\ref{cor:crtbp-verification}.
For $N=4$, \eqref{eq:Delta4} with $-d_{00}=g\delta$ gives \eqref{eq:line}, and on the axis $I_1=0$ one has $\ell_4(0,I_2)=d_{00}+d_{01}(0)I_2<0$.
\end{proof}

In the CRTBP the trivial branch $\eta=0$ consists of the planar Lyapunov family on the axis $I_2=0$ and the Lissajous family in the positive quadrant $I_1,I_2>0$.
Along the bifurcating branches, Proposition~\ref{prop:axis} turns a constraint root on the axis into an approximate halo orbit, which Remark~\ref{rem:axis-straightening} parameterizes with constant frequency.
Proposition~\ref{prop:invariance} turns a constraint root in the positive quadrant into an approximate quasi-halo torus, whose internal vector field \eqref{eq:innerflow} is angle dependent.

\subsection{Reconstruction in the CRTBP}
\label{ssec:fullreconstruction}

We now transport the approximate halo orbits and quasi-halo tori of Section~\ref{ssec:d4} to the CRTBP, where the full residual combines the center-manifold residual and the reduced residual.

\begin{corollary}[reconstruction of the quasi-halo tori]
\label{cor:fullreconstruction}
Assume the setting of Lemma~\ref{lem:cm} and Proposition~\ref{prop:invariance}, with $r$ as in Proposition~\ref{prop:invariance}, and let $K_N(\mathbb T^2)$ lie in the ball $B_{r_{\rm red}}^4$ with $r_{\rm red}\le r_{{\rm cm},N}/2$.
Define
\begin{equation}\label{eq:fullKN}
 K_N^{\rm full}:=E_N\circ K_N:
 \mathbb T^2\longrightarrow\mathbb R^6.
\end{equation}
Then the full residual of $K_N^{\rm full}$ decomposes pointwise as
\begin{align}
 &X_{\mathcal H}(K_N^{\rm full})
 -DK_N^{\rm full}v_N \notag\\
 &\quad=e_{{\rm cm},N}\circ K_N
 +(DE_N\circ K_N)
   \big(X_{H_{\le N}}(K_N)-DK_Nv_N\big).
 \label{eq:fullcovariant}
\end{align}
Consequently, with $M_{{\rm cm},N}:=\|DE_N\|_{C^0(B_{r_{\rm red}}^4)}$,
\begin{equation}\label{eq:fulldefect}
 \|X_{\mathcal H}(K_N^{\rm full})
 -DK_N^{\rm full}v_N\|_{C^0}
 \le
 C_{{\rm cm},N}r_{\rm red}^N
 +M_{{\rm cm},N}\widetilde C_Nr^N.
\end{equation}
Moreover, the two tori obtained for the roots $\pm\eta$ are exchanged by the physical vertical reflection and have the same internal dynamics up to the shift $\varphi_2\mapsto\varphi_2+\pi$.
They are the northern and southern finite-order quasi-halo tori.
\end{corollary}

\begin{proof}
Inserting and subtracting $(DE_N\circ K_N)X_{H_{\le N}}(K_N)$, using \eqref{eq:cmdefect}, and applying the chain rule to $E_N\circ K_N$ gives \eqref{eq:fullcovariant}.
The bound \eqref{eq:fulldefect} then follows from \eqref{eq:cmdefectbound} and \eqref{eq:invariance-defect}.
For the reflection statement, Proposition~\ref{prop:invariance} shows that $\Sop$ exchanges the two tori in the reduced coordinates and that their internal dynamics agree up to this shift.
Lemma~\ref{lem:cm} shows that $E_N$ intertwines $\Sop$ with the vertical reflection, so the exchange carries over to the CRTBP.
\end{proof}

If an angular straightening \eqref{eq:straightening} of $v_N$ is obtained independently, the full residual satisfies the same estimate, as in Remark~\ref{prop:rotation}.

\begin{corollary}[reconstruction of the halo orbits]
\label{cor:axisfull}
Assume the setting of Lemma~\ref{lem:cm} and Proposition~\ref{prop:axis}, with $r$ as in Proposition~\ref{prop:axis}, and let $K_N^{\rm ax}(\mathbb T^1)$ lie in the ball $B_{r_{\rm red}}^4$ with $r_{\rm red}\le r_{{\rm cm},N}/2$.
Then $K_N^{\rm ax,full}:=E_N\circ K_N^{\rm ax}$ satisfies
\begin{equation}\label{eq:axisfulldefect}
 \big\|X_{\mathcal H}(K_N^{\rm ax,full})
 -DK_N^{\rm ax,full}\,v_N^{\rm ax}\big\|_{C^0}
 \le
 C_{{\rm cm},N}r_{\rm red}^N
 +M_{{\rm cm},N}\widetilde C_Nr^N .
\end{equation}
Moreover, the two circles obtained for the roots $\pm\eta$ are exchanged by the physical vertical reflection and have identical tangential speed.
They are the northern and southern finite-order halo orbits.
\end{corollary}

\begin{proof}
The proof of Corollary~\ref{cor:fullreconstruction} applies with $K_N$, $v_N$, \eqref{eq:invariance-defect}, and Proposition~\ref{prop:invariance} replaced by $K_N^{\rm ax}$, $v_N^{\rm ax}$, \eqref{eq:axis-defect}, and Proposition~\ref{prop:axis}(iv).
\end{proof}

\subsection{A numerical example: the Sun--Earth
\texorpdfstring{$L_1$}{L1}}
\label{ssec:example}

For the Sun--Earth $L_1$, with $\mu=3.040423398444176\times10^{-6}$, the linear constants are listed in Table~\ref{tab:linear}.
The computations in this subsection use a C++ implementation of the recursion in double-precision arithmetic.
For $N\le20$ the condition of Section~\ref{ssec:verification} gives the admissible interval $|\eta|\le0.32$, which contains all roots used below.
\begin{table}[htbp]
\centering
\caption{Linear constants for the Sun--Earth $L_1$.}
\label{tab:linear}
\footnotesize % scriptsize
\setlength{\tabcolsep}{5pt}
\begin{tabular}{ccccccccccc}
\toprule
$c_2$ & $\omega_1$ & $\omega_2$ & $\delta$ & $g$ & $d_{00}$ & $\xi$ & $s_\omega$ & $\mx$ & $\mz$ & $\chixz$ \\
\midrule
$4.06107$ & $2.08645$ & $2.01521$ & $0.07124$ & $4.10166$ & $-0.29221$ & $-3.2293$ & $17.386$ & $0.33917$ & $0.99622$ & $0.16894$ \\
\bottomrule
\end{tabular}
\end{table}

The computed fourth-order counterterm coefficients are
\[
d_{10}(0)=1.771746284,\qquad d_{01}(0)=-1.592416172,
\]
so that
\[
I_1^{\min}=\frac{g\delta}{d_{10}(0)}=0.1649301946 .
\]
Since $d_{01}(0)<0$, the critical line does not meet the vertical axis $I_1=0$.
In terms of the planar amplitude $\alpha=\mx\sqrt{I_1}$ the threshold is $\alpha_{\min}=0.13774$, close to the halo thresholds obtained from resonant normal forms and from Lindstedt--Poincar\'e series \cite{LinChiba,Luo}.

At the threshold point $(I_1^{\min},0)$,
\[
\zeta_4(I_1^{\min},0)=-9.3484708069 .
\]
Corollary~\ref{thm:line-crtbp} therefore applies.
Since $\zeta_4<0$, \eqref{eq:localu} shows that the local pair of constraint roots $\eta=\pm\sqrt{u(I)}$ emerges on the side $\ell_4(I)>0$, namely $I_1>I_1^{\min}$ on the axis $I_2=0$ and, in the positive quadrant,
\[
0<I_2< \frac{d_{10}(0)}{|d_{01}(0)|}\,\big(I_1-I_1^{\min}\big).
\]

\begin{figure}[tbp] % [!t]
	\centering
	\begin{subfigure}[t]{0.48\textwidth}
		\centering
		\includegraphics[width=\linewidth]{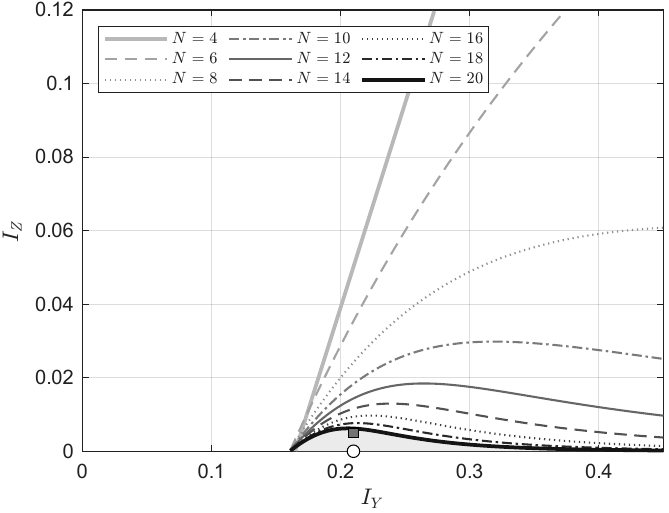}
		\caption{}
		\label{fig:1a}
	\end{subfigure}
	\hfill
	\begin{subfigure}[t]{0.49\textwidth}
		\centering
		\includegraphics[width=\linewidth]{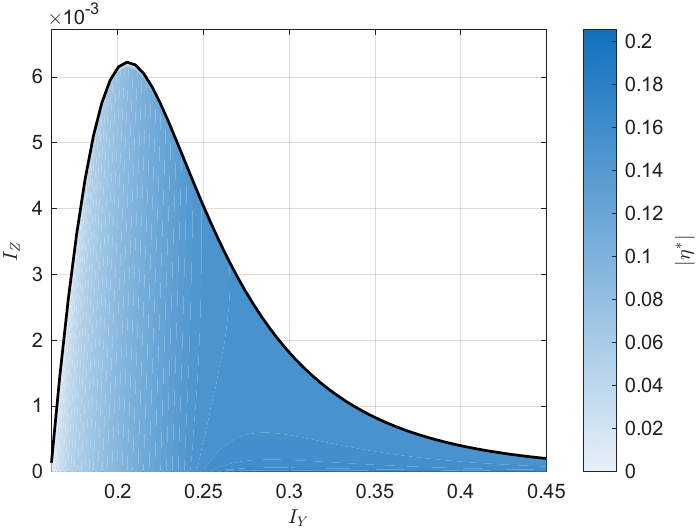}
		\caption{}
		\label{fig:1b}
	\end{subfigure}
	\caption{The bifurcation equation for the Sun--Earth $L_1$.
	(a)~Critical curves $\ell_N(I)=0$ for $N=4,6,\ldots,20$, with the two action points of Figure~\ref{fig:orbits}.
	(b)~Smallest positive root $|\eta_{20}^*|$ in the region where the $N=20$ bifurcating branches exist.
	The shading in (a) marks the domain of (b).}
	\label{fig:bifurcation}
\end{figure}

\begin{figure}[htbp]
	\centering
	\begin{subfigure}[t]{0.48\textwidth}
		\centering
		\includegraphics[width=\linewidth]{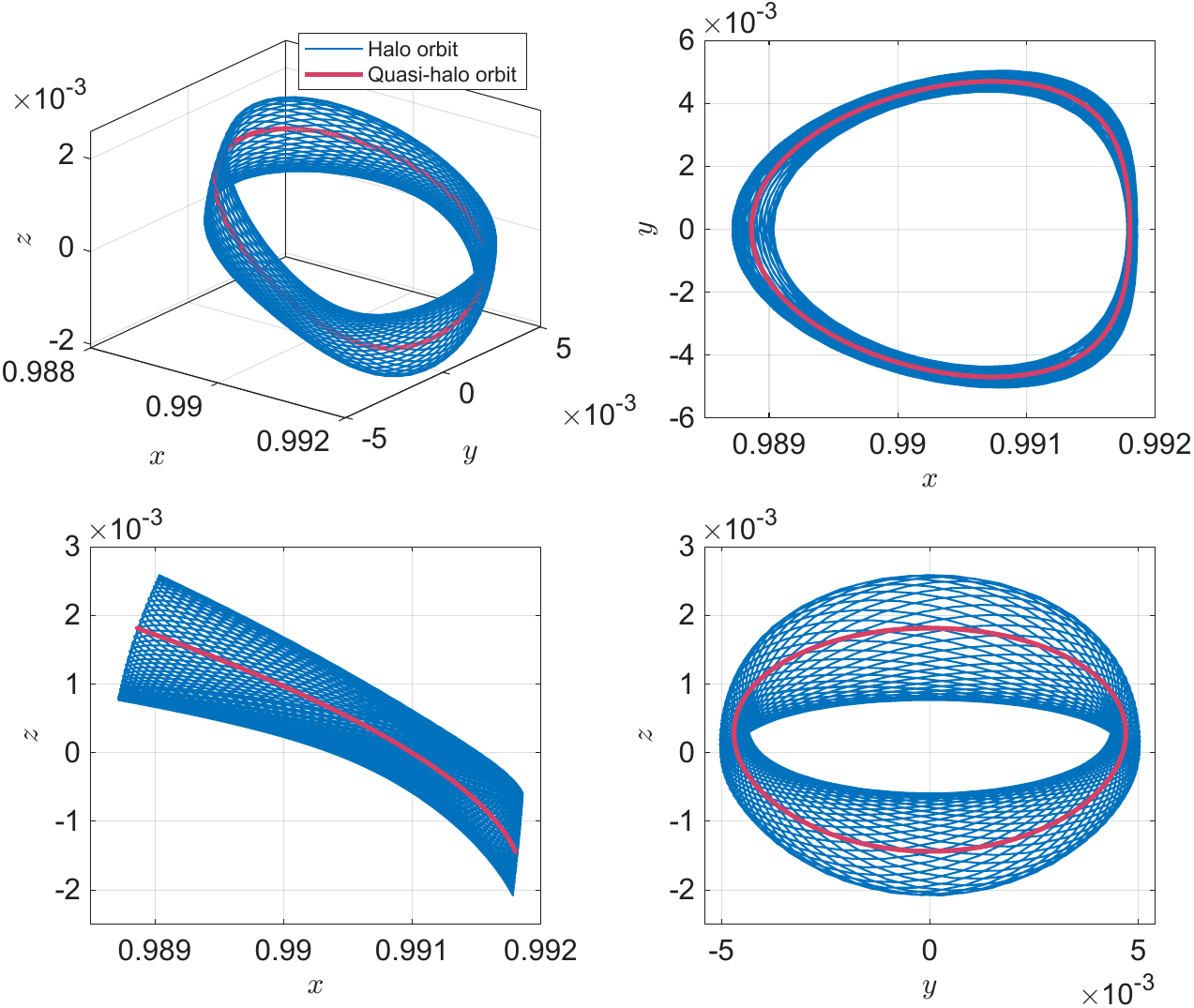}
		\caption{}
		\label{fig:2a}
	\end{subfigure}
	\hfill
	\begin{subfigure}[t]{0.49\textwidth}
		\centering
		\includegraphics[width=\linewidth]{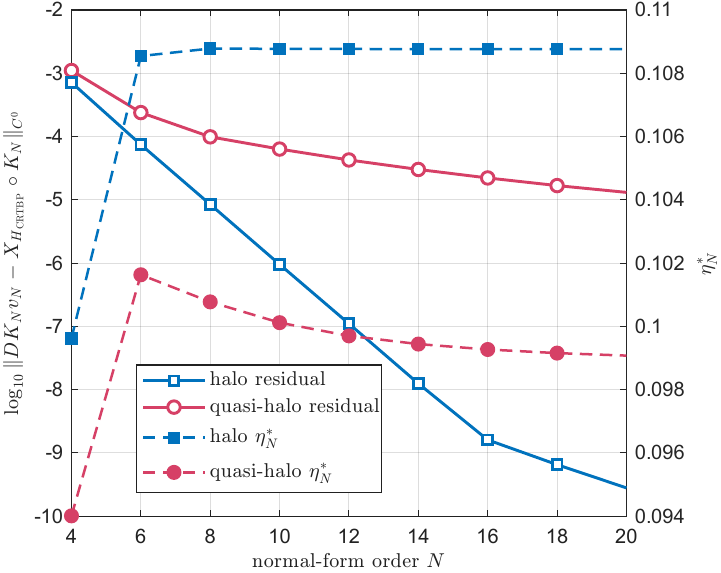}
		\caption{}
		\label{fig:2b}
	\end{subfigure}
	\caption{Reconstructed halo and quasi-halo orbits for the Sun--Earth $L_1$.
	(a)~The halo orbit at $(I_1,I_2)=(0.21,0)$ and the quasi-halo orbit at $(0.21,0.005)$, reconstructed at $N=10$ with $\eta_{10}^*=0.1087726$ and $0.1001184$, respectively.
	(b)~Full residual (left axis) and root $\eta_N^*$ (right axis) for $N=4,6,\ldots,20$.}
	\label{fig:orbits}
\end{figure}

At the two representative action points $(0.21,0)$ and $(0.21,0.005)$, which lie on the side $\ell_4(I)>0$, the fourth-order bifurcation equation has the roots $\eta_4^*=\pm0.0996$ and $\pm0.0940$.
The roots at the axis point $(0.21,0)$ illustrate Proposition~\ref{prop:axis} and Corollary~\ref{cor:axisfull} numerically, and those at the point $(0.21,0.005)$ illustrate Proposition~\ref{prop:invariance} and Corollary~\ref{cor:fullreconstruction}.

Figure~\ref{fig:bifurcation} shows the critical curves $\ell_N(I)=0$ of Theorem~\ref{thm:line} for the even orders $N=4,\ldots,20$, together with the distribution of the smallest positive root $\eta_{20}^*$ over the region where the $N=20$ bifurcating branches exist.
Because new counterterm coefficients $d_{ab}$ enter at each higher even order, the reconstruction at each even order $N=4,6,\ldots,20$ in Figure~\ref{fig:orbits} uses a root $\eta_N^*$ of $\Delta^{(N)}(I^*;\eta)=0$, continued from the preceding order.
At $N=20$ the bifurcation equation at the point $(0.21,0.005)$ has a second positive root, but the plotted value remains on the continued bifurcating branch.

For the reconstruction diagnostic we report the pointwise full residual.
The deviation from a numerically integrated trajectory is less informative, because the hyperbolic eigenvalue $\lambda\approx2.53$ amplifies integration errors by a factor of order $e^{\lambda T}$ over one period $T$.
In Figure~\ref{fig:orbits}(b) the residual of the halo orbit decreases from $7.2\times10^{-4}$ at $N=4$ to $2.8\times10^{-10}$ at $N=20$, and its root $\eta_N^*$ settles from $N=10$ on.
For the quasi-halo orbit the residual decreases more slowly, from $1.1\times10^{-3}$ to $1.3\times10^{-5}$, and $\eta_N^*$ still varies in the fourth significant digit at $N=20$.

% \FloatBarrier

% =====================================================================
\section{Discussion and outlook}\label{sec:discussion}
% =====================================================================

This paper constructs the coupling normal form, in which the first slow obstruction along a selected coupling direction is absorbed into an action-dependent counterterm.
Triangularity and symmetry propagation make the coupled construction uniquely solvable, and Theorem~\ref{thm:main} yields an action normal form with the counterterm factorization.
The resulting recovery constraint separates the trivial branch, described by the Birkhoff normal form, from the nontrivial branch defined by the bifurcation equation.
Theorem~\ref{thm:line} gives, at every finite order, the critical curves from which the bifurcating branches emerge.
Along these branches the finite-order construction provides approximately invariant tori and axis circles with explicit residual estimates, and in the CRTBP these become approximate quasi-halo tori and halo orbits.

The construction gives a Hamiltonian formulation of the coupling-induced bifurcation mechanism of \cite{LinChiba,LLC}, in which the recovery constraint is the counterpart of $\check\eta\,\check\Delta=0$ (Remark~\ref{rem:LCdict}).
Their series solutions locate the families and the bifurcation, and the coupling normal form adds analytic normalizing maps and residual bounds, which rigorous existence arguments require.

The formal theorem does not imply convergence of the normalizing map, and such convergence would amount to local analytic complete integrability of $H+\eta\Kc$.
The main open problem is the persistence of the approximately invariant objects, a prerequisite for a dynamical bifurcation.
A posteriori KAM and parameterization methods start from residual bounds of the type \eqref{eq:fulldefect} and \eqref{eq:axisfulldefect} \cite{DLLGJV,HaroBook,FontichDLLS}.
Rigorous implementations require validated estimates \cite{FiguerasHaroLuque}, as in \cite{CapinskiZgliczynski} for the restricted three-body problem.
For the axis circles such an argument can start directly (Remark~\ref{rem:axis-straightening}), whereas for the tori the internal vector field \eqref{eq:innerflow} must first be straightened, which is a small-divisor problem (Remark~\ref{prop:rotation}).
A dynamical bifurcation requires in addition a Floquet or variational degeneracy of the family on the trivial branch at the critical curve (Remark~\ref{rem:critical}).

\section*{Acknowledgments}
The author is grateful to Prof. Hayato Chiba for useful suggestions and comments.
The author was supported by the Sustainability Open Knowledge-Action Program by Connecting Multi-stakeholder Projects from Tohoku University and JST Moonshot, Grant Number JPMJMS2023, Japan.

\appendix

% =====================================================================
\section{Linear normalization and center manifold reduction}\label{app:linear}
% =====================================================================

This appendix proves Lemmas~\ref{lem:linlayer} and~\ref{lem:cm}, which follow the standard preprocessing of \cite{GJMS,JM}.

\subsection{Linear symplectic normalization}\label{app:linnorm}

This subsection proves Lemma~\ref{lem:linlayer}.
The vertical oscillator is normalized by $z=\omega_2^{-1/2}\bQ_2$ and $p_z=\omega_2^{1/2}\bP_2$.

The real invariant subspace of the planar block for the eigenvalues $\pm i\omega_1$ is spanned by
\[
b_1=(1,0,0,1+\xi\omega_1), \qquad b_2=(0,\xi,-(\omega_1+\xi),0)
\]
in the coordinates $(x,y,p_x,p_y)$, and the canonical symplectic form $dx\wedge dp_x+dy\wedge dp_y$ takes the value $-s_\omega$ on $(b_1,b_2)$.
Since $\xi<0$,
\[
s_\omega=\omega_1(1+\xi^2)-2|\xi| \ge 2|\xi|(\omega_1-1)>0,
\]
where $1+\xi^2\ge2|\xi|$ and $\omega_1>1$.
Hence $(b_1,-b_2)/\sqrt{s_\omega}$ is a symplectic basis of this subspace, and the planar part of $L$ is $(x,y,p_x,p_y)=(\bQ_1b_1-\bP_1b_2)/\sqrt{s_\omega}$.

The hyperbolic pair $(\bQ_h,\bP_h)$ is normalized as in \cite{JM}, by eigenvectors for the eigenvalues $\pm\lambda$ divided by $s_\lambda$ with $s_\lambda^2=2\lambda[(4+3c_2)\lambda^{2}+4+5c_2-6c_2^{2}]$.
The first components of the planar and vertical parts give \eqref{eq:xzrows0} with the constants \eqref{eq:mxmz}.

The reversor fixes $b_1$ and $z$ and negates $b_2$ and $p_z$, so it preserves the center subspace and acts on it as $\Rop$.
The vertical reflection changes only the signs of $z$ and $p_z$, so it preserves the center subspace and acts on it as $\Sop$.

\subsection{Center manifold reduction}\label{app:cm}

This subsection proves Lemma~\ref{lem:cm}.
For a monomial with unequal hyperbolic exponents the divisor of the center manifold reduction has real part equal to a nonzero integer multiple of $\lambda$, so its modulus is at least $\lambda$.
Successive Lie transforms therefore eliminate all such monomials without small-divisor conditions on the center frequencies.
Every remaining monomial contains $\bQ_h$ and $\bP_h$ with equal exponents, so the hyperbolic components of the vector field of $\mathcal H^{\rm cm}_{\le N}$ vanish on $\mathcal C$.
Hence $\mathcal C$ is invariant, and $X_{\mathcal H^{\rm cm}_{\le N}}\circ\iota=D\iota\,X_{H_{\le N}}$.

Let
\[
R^{\rm cm}_{N+1}:=(\mathcal H\circ L)\circ\Phi_{\rm cm}^{(N)}-\mathcal H^{\rm cm}_{\le N}.
\]
Lemma~\ref{lem:finite-lie} with $d=3$ and $F=\mathcal H\circ L$ gives the radius $r_{{\rm cm},N}$ and $\|X_{R^{\rm cm}_{N+1}}\|_r\le Cr^N$.
Since $L\circ\Phi_{\rm cm}^{(N)}$ is symplectic, the definition of $R^{\rm cm}_{N+1}$ and the identity for $X_{\mathcal H^{\rm cm}_{\le N}}\circ\iota$ give
\[
e_{{\rm cm},N}=\big(D(L\circ\Phi_{\rm cm}^{(N)})\circ\iota\big)\,X_{R^{\rm cm}_{N+1}}\circ\iota,
\]
and the uniform bound on this derivative proves \eqref{eq:cmdefectbound}.

Taking each generator without components on monomials with equal hyperbolic exponents makes it unique.
Each $\Gamma_n$ is then determined by $\mathcal H$ and the lower-degree generators, so it does not depend on $N$, and the $H_{\le N}$ are the truncations of one formal reduced Hamiltonian $H$.
The reversor exchanges the two hyperbolic eigendirections and the vertical reflection fixes them, so both involutions preserve the class of monomials with equal hyperbolic exponents.
The inherited symmetries then follow as in the proof of Proposition~\ref{prop:equivariant}(i), (ii), and~(v).

\end{document}